\documentclass{article}
\usepackage{amsfonts}
\usepackage{amssymb}
\newcommand{\iint}{\int\!\!\int}

\newcommand{\te}{Teich\-m\"{u}l\-ler}
\newcommand{\tes}{Teich\-m\"{u}l\-ler's}

\newenvironment{definition}{\bigskip \noindent {\sc
Definition.}}{\vskip
.2in}

\newtheorem{lemma}{Lemma}
\newtheorem{thm}{Theorem}

\newenvironment{proof}{{\sc Proof.}}{$\;\square$ \vskip .2in}
\begin{document}
\title{Universal \te\ Space}
\begin{author}{Frederick P. Gardiner and William J. Harvey}
\end{author}
\begin{date}{\ }
\end{date}
\maketitle

  %
  \renewcommand{\thefootnote}{\arabic{footnote}}%

\def\esssup{\mathop{\rm esssup}}

\centerline{\bf Abstract}

\medskip

We present an outline of the theory of universal \te\ space,
viewed as part of the theory of $QS,$ the space of quasisymmetric
homeomorphisms of a circle.   Although elements of $QS$ act in one
dimension, most results about $QS$ depend on a
two-dimensional proof.
$QS$ has a manifold structure modelled on a Banach space, and
after factorization by $PSL(2,{\mathbb R})$ it becomes a complex
manifold. In applications, $QS$ is seen to contain
many deformation spaces for dynamical systems acting in one, two
and three dimensions; it also contains deformation spaces of
every hyperbolic Riemann surface,
and in this naive sense it is universal. The deformation spaces are
complex
submanifolds and often have certain universal properties
themselves, but
those properties are not the object of this paper. Instead we focus
on the analytic foundations of the theory necessary for
applications to dynamical systems and rigidity.
\par
We divide the paper into two parts.  The first part concerns
the real theory of $QS$ and results that can be stated purely in
real terms; the basic properties are given mostly without
proof, except in certain cases when an easy real-variable
proof is available.
The second part of the paper brings in the complex analysis and
promotes
the view that properties of quasisymmetric maps are most easily
understood by
consideration of their possible two-dimensional quasiconformal
extensions.

\medskip

\begin{center}{\bf Table of Contents}
\end{center}

\noindent {\bf Introduction}

\medskip

\noindent {\bf 1.} {\bf Real Analysis}

\medskip

\indent 1.1 \hspace{.1in} Quasisymmetry

\indent 1.2 \hspace{.1in} The Quasisymmetric Topology

\indent 1.3 \hspace{.1in} The Symmetric Subgroup

\indent 1.4 \hspace{.1in} Dynamical Systems and Deformations

\indent 1.5 \hspace{.1in} The Tangent Spaces to $QS$ and $S$

\indent 1.6 \hspace{.1in} The Hilbert Transform and Almost Complex
Structure

\indent 1.7 \hspace{.1in} Scales and Trigonometric Approximation

\indent 1.8 \hspace{.1in} Automorphisms of \te\ Space

\medskip

\noindent {\bf 2.} {\bf Complex Analysis}

\medskip

\indent 2.1 \hspace{.1in} Quasiconformal Extensions

\indent 2.2 \hspace{.1in} \tes\ Metric

\indent 2.3 \hspace{.1in} Quadratic Differentials

\indent 2.4 \hspace{.1in} Reich-Strebel Inequalities

\indent 2.5 \hspace{.1in} The Tangent Spaces Revisited

\indent 2.6 \hspace{.1in} The Hilbert Transform and Almost Complex
Structure

\indent 2.7 \hspace{.1in} Complex Structures on Quasi-Fuchsian
Groups

\indent 2.8 \hspace{.1in} Automorphisms are Geometric

\vspace{.5in}

\noindent {\bf Introduction}

\bigskip
The origins of this theory lie in the study of deformations of
complex structure in
spaces of real dimension 2 and the moduli problem for Riemann
surfaces. It seems appropriate, therefore, to begin with a brief
sketch of
how the notion of a \te\ space first
arose, within this problem of variation of complex structure on a
topologically fixed compact Riemann surface.
For brevity we shall restrict attention to {\em hyperbolic}
Riemann
surfaces, which have as universal covering space the
unit disc; the terminology refers to the fact that, via
projection form Poincar\'e's Riemannian metric on the disc, all the
surfaces
are endowed with a structure of hyperbolic geometry.
\par

The definition of \te\ space stands out clearly as a key
stage in the struggle to justify, and to make precise, the
famous assertion of Riemann [{\it Theorie der Abel'schen
Functionen}, Crelle J., t. 54, (1857)] that the number of
(complex) parameters (or `moduli') needed to describe all surfaces
of genus $g\geq 2$ up to conformal equivalence is $3g-3$. After
preliminary work over many years by a substantial number of
eminent mathematicians, including F. Schottky, A. Hurwitz, F.
Klein, R. Fricke and H. Poincar\'e, the crucial new idea  was
introduced by O. \te\ around 1938, \cite{Teich1}, \cite{Teich2},
following earlier work of H.
Gr\"otzsch.  One specifies a {\it topological marking} of the base
surface and then considers all homeomorphisms to a target Riemann
surface which have the property that they distort the conformal
structure near each point by at most a bounded amount, using a
precise analytic measure of the distortion to be defined below.
Gr\"otzsch (see \cite{Gr1},\cite{Gr2}) had used
this approach to resolve similar problems in estimating distortion
for smooth mappings between plane domains; the term {\it
quasiconformal} was coined by L.V. Ahlfors around 1930 for the
class of
homeomorphisms to be employed. The method was strengthened,
generalised and applied to the case of closed Riemann surfaces with
striking effect  by \te , as we indicate below.
\par

A fundamental relationship exists between the quasiconformal
homeomorphisms of the hyperbolic disc and the induced boundary
maps of the circle, and this lies at the heart of the viewpoint on
Teichm\"uller theory to be presented here: for a general Riemann
surface, one must consider not
only deformations of complex structure in the interior  but also
the  ways in which the conformal structure may
change relative to the boundary. It turns out that both aspects
are best studied on the universal covering surface, the unit disc
$\Delta=\{|z|<1\}$: quasiconformal mappings of the disc extend to
homeomorphisms of the closed disc and many (but not all) of the
properties of a quasiconformal homeomorphism can be expressed
solely in terms of the boundary homeomorphism of the circle
induced by it.
\par

Let $QS$ be the space of sense-preserving, quasisymmetric
self-maps
of the unit circle; such maps are precisely
those occurring as the boundary value of some quasiconformal self
map of the disc $\Delta.$  A map
$H:\Delta\to \Delta$ is called {\it quasiconformal} (sometimes
abbreviated
to q-c) if $K(H)<\infty$, where $K(H)$ is
the essential supremum,  for $z\in \Delta$ , of
the local dilatations $K_z(H)$, and the
{\it local dilatation} $K_{z}(H)$ at $z$ is defined as
$$\limsup_{\epsilon
\rightarrow 0} \frac{
\max_{\theta}\{|H(z+\epsilon e^{i\theta})-H(z)|\}}{
\min_{\theta}\{|H(z+\epsilon e^{i\theta})-H(z)|\}},$$
which may be interpreted as the upper bound of local distortion
as measured on circles centered at $z$; compare with the
definition (\ref{Dil}) in \$2.  The set of all possible
quasiconformal extensions $H: \Delta \to \Delta$ of a given
quasisymmetric map $h$ may be regarded as the {\it mapping class}
of $h$ in the disc, and a mapping
$H_0$ is called {\it extremal} for its class if $K(H_0) \leq K(H)$
for every extension $H$ of $h.$
This notion of extremality for a mapping (within a homotopy class
of quasiconformal maps between two plane regions)
was also introduced by Gr\"otzsch (op.cit.), but it was
 \te\ who recognized the significance of
extremal maps in the study of deformations of complex structures.
He
applied them decisively in \cite{Teich1}, \cite{Teich2}, to
establish a
measure of distance between two marked surfaces: here the upper
bound for the
local distortion of the mapping over the base surface is to be
minimised.
\par
A base (hyperbolic) Riemann surface is given as the quotient space
$X_0= \Delta/\Gamma$
of the unit disc
under the action of a Fuchsian group $\Gamma,$ which is
by definition a discrete group of M\"obius transformations
which are conformal automorphisms of the disc;
if the group is torsion-free, then topologically,
$\Gamma\cong \pi_1(X_0)$ represents the group of deck transformations of the
universal covering
projection from $\Delta$ to $X_0$.
Suppose now that we are given a quasisymmetric
map $h$ of the circle with the property that
the conjugate group $\Gamma_1=h \circ \Gamma \circ h^{-1}$ is also
Fuchsian: the two orbit spaces $\Delta/\Gamma$ and
$\Delta/\Gamma_1$ can be
viewed as the same topological surface but with
different complex structures. The
{\em mapping class} of $h$ for the group $\Gamma$ is the subset of
its
mapping class (the set of {\it all} q-c self-mappings of the disc
extending the map $h$) consisting of
those q-c extensions $H$ of $h$  with the property that
every element $H \circ \gamma \circ H^{-1}$  of $\Gamma_1$ acts as
a M\"obius
transformation of the disc $\Delta$ to itself.
For a Fuchsian group $\Gamma$ that covers a compact Riemann
surface, \tes\ theorem establishes a profound link
between an extremal representative $H$ for a given class and a
holomorphic quadratic differential for $\Gamma$ - a complete proof
is given in \cite{Bers1}.
As a consequence, one may infer that the
space of marked deformations
of the compact surface $\Delta/\Gamma$ is a complete metric space
homeomorphic to an $n=(6g-6)$-dimensional real cell.
The metric is called \tes\ metric and the distance between the
base surface $X_0=\Delta/\Gamma$ and the marked surface
$X_1=\Delta/\Gamma_1$, with $\Gamma_1=H_0 \circ \Gamma \circ
H_0^{-1}$,
is $\log K(H_0),$ where
$H_0$
is extremal in its class.
This type of
extremal mapping is a feature of continuing interest, partly
because of the connection with Thurston's theory of measured
laminations on hyperbolic surfaces, \cite{Bers4},\cite{EG1},
\cite{HM},
\cite{G}, \cite{Harvey2},\cite{S2}.
\par

The final ingredient, which makes it possible to construct
these holomorphic parameter spaces
for all types of Riemann surface, is the relationship
between the quasiconformal property and the solutions of a certain
partial
differential equation.
By a fundamental observation of Lipman Bers (see \cite{AB}), if $H$
is a
quasiconformal self-map of the disc, it
satisfies the Beltrami equation
\begin{equation}
\label{Beltrami equation}
H_{\overline{z}}(z) = \mu(z) H_z.
\end{equation}
where $\mu$, with $||\mu||_{\infty}<1$, is a measurable
complex-valued function
on the disc, which represents the complex dilatation at each point
of $\Delta$;
$\mu$ is often called the {\em Beltrami coefficient} of $H$.
Conversely,
by virtue of solvability properties of this equation, $\mu$
determines $H$  uniquely
up to postcomposition by a M\"obius transformation.
By using the analytic dependence of $H=H_{\mu}$ on its Beltrami
coefficient
$\mu$, and deploying a construction known as the Bers embedding,
each
$T(\Gamma)$ is embedded as a closed subspace of the
complex Banach space $B$ of univalent functions on $\Delta$ which
have quasiconformal
extensions to the sphere --  more
details are given in
section 2.5. It then follows
that $T(\Gamma)$ has a natural structure of complex manifold for
{\it any} Fuchsian group $\Gamma$. Furthermore, each inclusion
$\Gamma'\subset \Gamma$ of Fuchsian groups induces a contravariant
inclusion of these \te\ spaces $T(\Gamma)\subset T(\Gamma')$, which
implies that  the Banach space $T(1)=B$,  which corresponds to the
trivial Fuchsian group $\Gamma'= 1= \langle  \hbox{Id}  \rangle$
is {\em universal} in the sense that it contains the \te\ spaces
of every hyperbolic Riemann surface $\Delta/\Gamma$.
\par
In the period after World War II,  the verification of \tes\ ideas
and the
subsequent rigorous development of the foundational complex
analytic deformation theory outlined above by L. V.~Ahlfors,
L.~Bers, H.E.~Rauch and their students
occupied more than 20 years.
The circumstances of \tes\ life and particularly his
political activities caused much
controversy and, coupled with the relative inaccessibility of his
publications, this perhaps contributed to some early reluctance
to pursue a theory based on his claims; for commentary on
mathematical life in Germany under the Third Reich, the reader
might consult
\cite{Teich4} and \cite{Tietz}.
Detailed
expositions of this foundational work on moduli are given in
\cite{G}, \cite{E2},  \cite{IT}, \cite{Nag} and
\cite{Lehto}.
\par
In a landmark study of the local complex analytic geometry of \te\
space, H.L.~Royden \cite{Royden} showed
that when
$T(\Gamma)$ is finite dimensional,
the complex structure of the space determines its \te\
metric.  In fact, he
proved that \tes\ metric coincides with the Kobayashi metric
\cite{Kobayashi}, which is defined purely in terms of the set of
all holomorphic maps from the unit disc into
$T(\Gamma).$
Royden also showed that every biholomorphic automorphism of
$T(\Gamma)$
is induced geometrically by an element of the mapping class group,
a result which extends to many infinite dimensional \te\ spaces;
we examine this important rigidity theorem more carefully
in sections $1.8$ and $2.8$.
\par

The case of compact Riemann surfaces and their deformation spaces
calls for techniques involving aspects of surface topology
and geometry which will not be considered in this article.
Instead, we present a  formulation  which focusses on the real
analytic foundations
of the theory, important for applications to real and complex
dynamical systems and
matters which relate to rigidity.
It was observed by S.P.~Kerckhoff (see for instance \cite{Wolpert})
and
later, independently, by S.~Nag and A.~Verjovsky \cite{NV}
that the almost complex structure
on each $T(\Gamma)$ corresponding to its complex structure is given
by
the Hilbert transform acting on the relevant space of vector fields
defined on the unit
circle.
This fact indicates that deep results concerning the complex
structure of \te\ space can be viewed purely as theorems of real
analysis.
With this principle in mind, we divide the exposition into two
parts.  The first part
concentrates on
the real theory of $QS$ and we present the theorems
in real terms as far as possible; the basic properties are stated
mostly without
proof, except in certain cases where an easy real-variable
proof is available.
The second part of the paper follows closely the
outline of the first but brings in the complex analysis: in our
view,
despite their very real nature,
properties of quasisymmetric maps are most easily
understood by
consideration of their possible two-dimensional quasiconformal
extensions.

\begin{section}{Real Analysis}

\begin{subsection}{Quasisymmetry}

A quasisymmetric map $h$ of an interval $I$ to an interval
$J$ is an increasing homeomorphism $h$ for which there exists a
constant $M$
such that
\begin{equation}
\label{dfn1}
\frac{1}{M} \leq \frac{h(x+t)-h(x)}{h(x)-h(x-t)} \leq M
\end{equation}
for every $x$ and $t>0$ with $x-t, x$ and $x+t$ in $I.$ It is not
hard to prove the quasisymmetric maps form a pseudo-group. That
is, if $h$ is quasisymmetric from $I$ to $J$ with constant $M,$
then $h^{-1}$ from $J$ to $I$ is quasisymmetric with constant
$M_1$ depending only on $M.$  Moreover, if $g$ is
quasisymmetric from $I_1$ to $I_2$ with constant $M_g$ and $h$ is
quasisymmetric from $I_2$ to $I_3$ with constant $M_h,$ then $h
\circ g$ is quasisymmetric from $I_1$ to $I_3$ with constant $M_{h
\circ g}$ depending
only on $M_g$ and $M_g.$ Also, $h$ is H\"older continuous with
H\"older exponent $\alpha$ depending only on $M.$ For purposes of
illustration, we prove this fact here.

\begin{lemma} A quasisymmetric map of an interval $I$ to an
interval $J$
satisfying condition {\rm (\ref{dfn1})} is H\"older continuous.
\end{lemma}
\begin{proof} (See Ahlfors \cite{A}, page 65 and 66).
Pre and postcomposition of $h$ by affine maps yields a map
$\tilde{h}=A
\circ h \circ B$ with the same
constant $M$ of quasisymmetry.
To show $h$ is H\"older continuous at a point $p$, it suffices to
show
$\tilde{h}$ is H\"older continuous at $0$ and to assume
the intervals $I$ and $J$ contain $[0,1],$
and that $\tilde{h}(0)=0$ and $\tilde{h}(1)=1.$
By plugging in $x-t=0,x=1/2,x+t=1$
inequality (\ref{dfn1}) yields
$$1/(M+1) \leq \tilde{h}(1/2) \leq M/(M+1).$$  Repeated
applications of
(\ref{dfn1}) with $x=1/2^n$ and $t=1/2^n$ yield
$$ 1/(M+1)^n \leq \tilde{h}(1/2^n) \leq (M/(M+1))^n.$$
Since $\tilde{h}$ is increasing, this implies that for $1/2^n \leq
x
\leq 1/2^{n-1},$ $$\tilde{h}(x) \leq \tilde{h}(1/2^{n-1}) \leq
(M/(M+1))^{n-1}= \left((M+1)/M\right)\left(M/(M+1)\right)^n.$$
But since $x \geq 1/2^n,$ the previous inequality implies
$$\tilde{h}(x) \leq ((M+1)/M)x^{\alpha}, {\rm \ where \ }
\alpha=\frac{\log((M+1)/M))}{\log 2}.$$
To finish the proof we note
that $h$ and $\tilde{h}$ have the same H\"older exponent
$\alpha.$\end{proof}

A
quasisymmetric homeomorphism of the unit circle
$S^1=\{e^{i\theta}: \theta {\rm \ real \ }\}$ is an orientation
preserving homeomorphism of the circle for which there exists a
constant $M$ such that
\begin{equation}
\label{dfn2}
\frac{1}{M} \leq \left|
\frac{f(e^{i(x+t)})-f(e^{ix})}{f(e^{ix})-f(e^{i(x-t)})}\right|
\leq M,
\end{equation}
for all $x$ and all $|t|<\frac{\pi}{2}.$
 Obviously the restriction of any M\"obius transformation
preserving the unit disc to the unit circle is quasisymmetric. Let
a finite number of smooth real-valued charts $\varphi_j$ that cover
the
circle be given and assume the maps from intervals to intervals
defined
by $\varphi_j \circ h \circ (\varphi_k)^{-1}$ are quasisymmetric
with constants $M_{jk}$ on the intervals where they are defined.
Then $h$ will be quasisymmetric with a constant $M$ depending on
the $M_{jk}$ and the coordinate charts $\varphi_j$ and $\varphi_k.$
Conversely,
suppose $h$ is quasisymmetric and $\varphi_j$ are a finite system
of
smooth charts whose domains of definition cover the circle.  Then
each
$\varphi_j \circ
f \circ (\varphi_k)^{-1}$ is quasisymmetric on the interval where
it is
defined.  Thus, if when expressed in terms of a
finite number
of smooth charts that cover the circle $h$ is quasisymmetric on the
intervals for some constant $M,$ then $h$ is quasisymmetric on the
circle for some possibly different constant $M_1.$

\end{subsection}
\begin{subsection}{The Quasisymmetric Topology}
Now we introduce a topology on the group of orientation preserving
homeomorphisms $h$ of a circle that satisfy inequality
(\ref{dfn2}) by specifying a neighborhood basis $V(\epsilon)$
of the identity.  By definition $h$ is in $V(\epsilon),
\epsilon>0,$ if
two conditions are satisfied:

\vspace{.1in}

\begin{tabular}{l}
(a)  $\sup\{ \left|h(e^{ix}) - e^{ix}\right|,
\left|h^{-1}(e^{ix})-e^{ix}\right|\}< \epsilon,$ and\\
\  \ \\
(b)
inequality  (\ref{dfn2}) is satisfied with
$M=1+\epsilon$.
\end{tabular}

\vspace{.2in}

\noindent This system of neighborhood has the following properties:

\vspace{.1in}

\begin{tabular}{l}
(i) $\bigcap_{n=1}^{\infty} V(1/n) = \{identity\},$ \\
\ \\
(ii) for every $\epsilon>0,$ there exists $\delta>0,$ such that
$V(\delta) \circ V(\delta) \subset V(\epsilon),$ and \\
\ \\
(iii) for every $\epsilon>0,$ there exists $\delta>0,$ such that
$\left( V(\delta)\right)^{-1} \subset V(\epsilon).$
\end{tabular}

\vspace{.1in}

This system of neighborhoods induces a right and a left topology
on $QS$ by right and left translation.  That is, $V \circ h $ is
a
right neighborhood of $h$ when $V$ is a
neighborhood of the identity.  These neighborhoods are precisely
those that
make right
translation maps $h \mapsto h \circ g$ continuous.   Similarly,
there is the system of left neighborhoods $h \circ V$ of the
$h,$ and these make left translation maps $h
\mapsto g \circ h$ continuous.
However, these properties constitute only
part of the structure necessary
to make $QS$ a topological group.
In the next section we examine this discrepancy in more detail.

\end{subsection}
\begin{subsection}{The Symmetric Subgroup}
There is a brief, relatively undeveloped, theory of groups
that are also Hausdorff topological spaces satisfying axioms (i),
(ii),
and (iii) above. More details may be found in  \cite{GS}.
We summarize this theory and its application to $QS$ in
this section.

\begin{definition} A {\em topological group} is
a group $G$ that is also a Hausdorff topological space and such
that the map $(f,g) \mapsto f \circ g^{-1}$ from $G \times G$ to
$G$ is continuous.
\end{definition}

It turns out that $QS$ is not a topological group because taking
inverses
is not continuous.
However, it does satisfy the
axioms for what we call a partial topological group.

\begin{definition}
A {\em partial topological group} is a
group with a Hausdorff system
of neighborhoods of the identity satisfying {\rm (i), (ii)} and
{\rm (iii)}
above.
\end{definition}
As we have seen in section 1.2, at a general point $h$ of the group
there are two neighborhood systems.
If $U$ runs through the
neighborhood system at the identity, then $h \circ U$ and $U \circ
h$
run through systems of left and right neighborhoods of $h,$
respectively.
The following theorem is proved in \cite{GS}.

\begin{thm}
\label{thm1}
The following conditions on a partial topological group are
equivalent:

\vspace{.05in}

 \indent i{\rm )} it is a topological group with the given
neighborhood system of the \\
\indent identity,

\vspace{.05in}

\indent ii{\rm )}  the left and the right neighborhood systems
agree at every point,

\vspace{.05in}

\indent iii{\rm )} the adjoint map $f \mapsto h \circ f \circ
h^{-1}$ is continuous at the identity \\
 \indent  for every $h$ in the group.
\end{thm}

In a general partial topological group the properties of Theorem
\ref{thm1} will not hold.  One of the two topologies in a
partial topological group will be left translation invariant and
the
other right translation invariant.  The inverse operation
interchanges
these two topologies.

One can consider those elements $h$ of a partial topological group
for
which the two neighborhood systems at $h$ agree, that is,
for which conjugation by $h$ maps the neighborhood system at the
identity isomorphically onto itself.  These elements form a closed
subgroup:
the two topologies agree on this subgroup and give it the structure
of a
topological group.  We call this subgroup the {\em
characteristic topological subgroup}.

If a subset of a partial topological group is invariant under the
inverse operation, then it is closed for one topology if, and only
if,
it
is closed for the other.  In particular, one may speak without
ambiguity
of a closed subgroup of a partial topological group.

The next result is elementary.

\begin{thm}
\label{thm2}
The characteristic topological subgroup of a partial topological
group
is a closed topological subgroup.
 \end{thm}

\begin{definition}
\label{dfnxvi2}
A quasisymmetric map $h$ has {\em vanishing ratio distortion}
if there is a function $\epsilon(t)$ with $\epsilon(t)$
converging to zero as $t$ converges to zero,
such that inequality {\rm (\ref{dfn2})} is satisfied with $M$
replaced by $1+\epsilon(t).$
\end{definition}

 It turns out that the characteristic topological subgroup of
$QS$ comprises precisely those homeomorphisms that have vanishing
ratio  distortion.
We shall call this subgroup the {\em symmetric subgroup} $S.$
A direct proof that $S$ is a topological group is elementary.
Here we prove only the following fact.

\begin{thm} $S$ is a closed subgroup of $QS.$
\end{thm}
\begin{proof}  We shall use the following notation.  $I$ and $J$
are contiguous intervals, $I=[a,b], J=[b,c],$ and $|I|=b-a$ is the
length of $I.$  Let a constant $C>1$ be given.  One first shows
that  if $I$ and $J$ are contiguous with $$1/C \leq
|I|/|J| \leq C$$ and if $g$ is sufficiently near the identity in
the quasisymmetric topology, then
$$ \frac{1}{1+\epsilon} \leq
\frac{|g(I)|}{|g(J)|}\cdot \frac{|J|}{|I|} < 1+\epsilon.$$
Assume $h_n$ is a
sequence of elements of $S$ converging in the right $QS$-topology
to
$h.$  This means that for sufficiently large $n,$
$$\frac{1}{1+\epsilon} \leq \frac{|h \circ h_n^{-1}(I)|}{|h \circ
h_n^{-1}(J)|} \cdot \frac{|J|}{|I|} < 1 + \epsilon.$$
Also, assume that for each $h_n$ there is a function
$\epsilon_n(t)$ approaching zero as $t$ approaches $0$ such that
for all contiguous intervals $K$ and $L$ with $1/C \leq |K|/|L|
\leq C,$
$$\frac{1}{1+\epsilon_n(|K|)} <
\frac{|h_n(K)|}{|h_n(L)|}\cdot \frac{|L|}{|K|}
<1+\epsilon_n(|K|).$$
Taking the product, we obtain
$$\frac{1}{(1+\epsilon)(1+\epsilon_n(|K|))} \leq
\frac{|h \circ h_n^{-1}(I)|}{|h \circ h_n^{-1}(J)|} \cdot
 \frac{|J|}{|I|}\cdot \frac{|h_n(K)|}{|h_n(L)|}\cdot
 \frac{|L|}{|K|}
 < (1 + \epsilon)(1+\epsilon_n(|K|).$$
 Since there is a uniform bound on the quasisymmetric norm
 of $h_n,$ we may assume $1/C \leq \frac{|h_n(K)|}{|h_n(L)|} \leq
C,$
and thus we may substitute in $I=h_n(K)$ and $J=h_n(L).$ We obtain
 $$\frac{1}{(1+\epsilon)(1+\epsilon_n(|K|))} \leq
\frac{|h(K)|}{|h(L)|} \cdot
 \frac{|L|}{|K|}
 < (1 + \epsilon)\left(1+\epsilon_n(|K|)\right).$$
For given $\delta>0,$ we can pick $n_0$ large enough so that
$\epsilon_{n_0}(|K|)<\epsilon$ whenever $|K|<\delta$ and $1/C
<|K|/|L|<C.$  Then
$$\frac{1}{(1+\epsilon)^2} <\frac{|h(K)|}{|h(L)|} \cdot
 \frac{|L|}{|K|} < (1+\epsilon)^2,$$
 and this implies $h$ has vanishing ratio distortion.
\end{proof}
\end{subsection}

\begin{subsection}{Dynamical Systems and Deformations}
By definition, {\em universal Teichm\"uller space} ${T}$ is
$QS$ factored by the close subgroup of M\"obius transformations
that
preserve the unit disc.  It is universal in the naive sense that
it contains the deformation spaces of nearly all one-dimensional
dynamical systems $F$ that act on the unit circle.
When we say this, we have in mind the following two types of
dynamical systems. $F$ is either a
 Fuchsian group
acting on the unit circle or a $C^2$-homeomorphism acting on the
unit circle.   In the second situation, it is often
useful to assume $F$ has irrational rotation number.  The
theory is already complicated if $F$ is a diffeomorphism and
becomes even more so if $F$ is allowed to have one critical point.
Under smooth changes of
coordinate, we may assume $F$ maps an interval on the real axis to
another interval on the real axis and maps the origin to a point
$c.$
We also assume that in suitable smooth coordinates $F$ takes the
form of
a power law: $$ F(x) = |x|^{\alpha}{\rm sign}(x) + c,$$
for some constant $\alpha>1$.

The deformation space $T(F)$ (sometimes called the
Teichm\"uller space of $F$), is defined to be the space of
equivalence classes of quasisymmetric maps $h \in QS$ such that $h
\circ F \circ h^{-1}$ is a dynamical system of the same type.  Two
maps $h_0$ and $h_1$ are equivalent  if there is a M\"obius
transformation $A$ such that $A \circ h_0=h_1.$
Thus, it is the set of quasisymmetric conjugacies to
dynamical systems of the same type factored by this equivalence
relation.

In the case that $F$ is a Fuchsian group with generators $\gamma_j$
 this means that, for each $j,$
the conjugate $h \circ \gamma_j \circ h^{-1}$ is also a M\"obius
transformation preserving the unit circle.  In the case $F$ is
a $C^2$ homeomorphism possibly with a power law, this means $h
\circ
F \circ h^{-1}$ is also a $C^{2}$-homeomorphism.
If $h$ is itself a M\"obius
transformation, we consider the dynamical systems generated
by $F$ and by \linebreak
$h \circ F \circ h^{-1}$ as not differing in
any essential way.  For this reason, we view $T(F)$ as a subspace
of
$QS {\rm \ mod \ } PSL(2,\mathbb R).$  That is, two elements $h_1$
and
$h_2$ are considered equivalent if there is a M\"obius
transformation $A$ such that $A \circ h_1 = h_2.$

It turns out that the factor space $T=QS {\rm \ mod \ }
PSL(2,\mathbb R)$
carries in a natural way the structure of a complex manifold,
as do the subspaces $T(F)$ for
many dynamical systems $F.$  Even the statement
that $T(F)$ is connected is
already significant, and knowledge of geometrical properties
of curves which join pairs of points in $T(F)$
can have dynamical consequences.

Here we explain why
conjugacies that allow
distortion of eigenvalues cannot be smooth: therefore, to obtain
interesting
conjugacies one must expand into the quasisymmetric realm.

\begin{lemma} Let $F_0$ and $F_1$ be
two discrete dynamical systems acting on the real axis,
generated by $x
\mapsto \gamma_0(x)=\lambda_0 x$ and by $x \mapsto \gamma_1(x)=
\lambda_1 x,$ respectively,
and assume $1<\lambda_0 <\lambda_1.$ Let
$h$ be a conjugacy, so that
$h \circ \gamma_0 \circ h^{-1} = \gamma_1.$  Then $h$ can be at
most
H\"older continuous with exponent $\alpha =
\log \lambda_0/\log \lambda_1.$
\end{lemma}
\begin{proof}
Because $h(\lambda_0^n x) = \lambda_1^n h(x),$ by plugging in $x=1$
and letting $n$ approach $-\infty$ and $\infty,$ one sees that
$h$ must fix $0$ and $\infty.$  By postcomposition of $h$ with
a real dilation, we may assume $h(1)=1$ and this implies
$h(\lambda_0^n)=\lambda_1^n.$  But any such map taking these
values for arbitrarily large negative values of $n$ cannot satisfy
an
inequality of the form $|h(x)| \leq C|x|^{\alpha}$ unless
$\alpha \leq \log \lambda_0/ \log \lambda_1.$
\end{proof}
\end{subsection}

\begin{subsection}{Tangent Spaces to $QS$ and $S$}
In this section we identify the circle with the
extended real line $\overline{\mathbb R}.$
By postcomposing with a M\"obius transformation
we may assume any homeomorphism
representing an element of $T$ fixes infinity.
Consider a smooth curve $h_s$ of
homeomorphisms in $QS$ parameterized by $s$ and passing through
the identity at $s=0.$
We may assume each homeomorphism $h_s$ fixes infinity.
By smooth we shall mean that
\begin{equation}
\label{infin}
 h_s(x) = x + s V(x) + o(s),
\end{equation}
where the distance measured in the quasisymmetric norm from the
identity to $h_s$ is less than or equal to
 a constant times $s.$
In particular, $$ \frac{1}{1+C s} \leq
\frac{h_s(x+t)-h_s(x)}{h_s(x)-h_s(x-t)} \leq 1 + C s.$$ By
substituting (\ref{infin}) into this formula we arrive at the
following condition on the continuous function $V:$
\begin{equation}
\label{bigZ}
 \left|V(x+t)-2V(x) +V(x-t) \right| = O(t).
 \end{equation}
If $h_s$ is a smooth curve in the symmetric subspace $S,$ then
\begin{equation}
\label{littleZ}
 \left|V(x+t)-2V(x) +V(x-t) \right| = o(t).
 \end{equation}
We will call (\ref{bigZ}) and (\ref{littleZ}), respectively,
the big and little Zygmund conditions.
Since $V$ is to be regarded as the tangent vector to the
one-parameter family
of homeomorphisms $h_s,$ $V(x)\frac{\partial }{\partial x}$ is a
vector field.

 If instead we consider the mappings $h_s$ as acting on the
unit circle $|z|=1$ then the condition that the vector field $W$
point in a
direction tangent to the unit circle is that
\begin{equation}
\label{trigtrans}
\tilde{W}(x)
=W(e^{ix})/ie^{ix}
\end{equation}
be real-valued.  The boundedness conditions on
$QS$ and $S$ correspond to the conditions
 that the continuous, periodic function $\tilde{W}$
 satisfy (\ref{bigZ}) and (\ref{littleZ}).
 We denote the spaces of continuous vector fields satisfying these
conditions
by $Z$ and $Z_0,$ respectively.

A simple example of a tangent vector in $Z_0$ is
generated by a curve of M\"obius transformations preserving the
unit circle and passing through the identity. Such a curve has
a tangent vector of the form
$$ W(z) \frac{\partial}{\partial z} =
(\alpha z^2 + \beta z + \gamma) \frac{\partial}{\partial z},$$
where $\alpha,$ $\beta,$ and $\gamma$ are constants which make
$W(z)$ real-valued along $z \overline{z}=1.$  We call such tangent
vectors
trivial.  Thus the quadratic
polynomials which satisfy this reality condition
 define the tangent vectors to trivial curves of
homeomorphisms.

We will show that any tangent vector satisfying the big
Zygmund condition is the tangent vector to a smooth curve in $QS$
passing through the identity, and correspondingly any tangent
vector satisfying the little Zygmund condition is the tangent
vector to a smooth curve in $S.$

\begin{definition} Let
${\cal Z}$ and ${\cal Z}_0$ be the spaces $Z$ and $Z_0$ factored
by the quadratic polynomials.
\end{definition}

  Eventually in section 2.3 we shall identify a
Banach space ${\cal A}$ such that the Banach dual of ${\cal Z}_0$
is
isomorphic to ${\cal A}$ and the Banach dual of ${\cal A}$ is
isomorphic to ${\cal Z}.$ In particular, $Z_0^{**} \cong Z.$

Let
$Q$ be any quadruple of points $a,b,c$ and $d$ arranged in
counter-clockwise order on the unit circle or in increasing order
on the real axis and define the cross ratio $cr(Q)$ by
 \begin{equation}
\label{cr}
cr(Q)=\frac{(d-c)(b-a)}{(c-b)(a-d)}.
\end{equation}
Recall that $cr(Q)$ is M\"obius invariant in the sense that
$cr(A(q))=cr(Q)$ for any M\"obius transformation $A.$
In consequence we may define a norm $|| \ ||_{cr}$
on vector fields which is M\"obius invariant
in the sense that
$$||W||_{cr} = ||\frac{W \circ A}{A'}||_{cr},$$
for every M\"obius transformation $A.$

Define $W[a,b,c,d]$ to be the alternating sum
$$\frac{W(d)-W(c)}{d-c}-
\frac{W(c)-W(b)}{c-b}+
\frac{W(b)-W(a)}{b-a}-
\frac{W(a)-W(d)}{a-d}.$$
For a given quadruple $Q$ the term $cr(Q) \rho(cr(Q)) W[a,b,c,d]$
measures the velocity of the cross-ratio (\ref{cr})
with respect to the Poincar\'e metric $\rho(z)|dz|$
 on the sphere punctured at $0, 1$
and $\infty$ when each of the points $a, b, c$ and $d$ move with
complex velocities $W(a), W(b), W(c)$ and $W(d),$ respectively.
The {\em infinitesimal cross-ratio norm} is defined
for the space ${\cal Z}$  by
\begin{equation}
\label{crnorm}
||W||_{cr}= \sup_Q \left|  cr(Q) \rho(cr(Q))
W[a,b,c,d]
\right|.
\end{equation}

Note that $||W||_{cr}=0$ if, and only if, $W$ is a quadratic
polynomial.
Furthermore,  if $Q$ has the form $Q=(-\infty, x-t,x, x+t)$ then
$cr(Q)=-1.$
If in addition we assume $|W(z)|=o(|z|^2),$
which is tantamount to the assumption that
$W(z)\frac{\partial}{\partial z}$ vanishes
at infinity,
then the alternating sum $W[a,b,c,d]$ is equal to
$$- \frac{W(x+t)-2W(x)+W(x-t)}{t}.$$

\end{subsection}
\begin{subsection}{The Hilbert Transform and Almost Complex
Structure}
If the vector field
$\tilde{V}(z)\frac{\partial}{\partial z}$ is
continuous and
real-valued on the circle, then the function $V(x) =
\tilde{V}(e^{ix})/ie^{ix}$
is continuous, real-valued and periodic on the real axis.
Define a function $W(x)$ by the formula

\begin{equation}
\label{Hilb} W(x) = \frac{1}{\pi} \lim_{\epsilon \rightarrow 0}
 \int_{|y-x| \geq \epsilon} V(y) \cot (\frac{y-x}{2}) \  dy ,
\end{equation}
where the integral is taken over values of $y$ for which $-\pi
\leq y \leq \pi$ and $|y-x \pm 2 \pi n| \geq \epsilon$ for all
integers $n.$ Transporting $W$ back to the unit circle by the
formula $$\tilde{W}(e^{ix})=W(x)ie^{ix},$$ one obtains a
complex-valued function $\tilde{W}$ defined on the circle for
which $$\tilde{W}(z)\frac{\partial}{\partial z}$$ is again
real-valued. By this process $\tilde{V}$ is transformed to $J
\tilde{V} = \tilde{W},$ another field of vectors on the
unit circle whose directions are tangent to the circle.

This rule defines an operator $J$  called the {\em Hilbert
transform};
it extends to a bounded operator for many different
smoothness
classes.  For example,
$$||W||_p \leq C_p ||V||_p$$
for $p \geq 2,$
where $$||V||^p_p = \frac{1}{2\pi}
\int_{-\pi}^{\pi}|V(x)|^p dx.$$
More important to us are the
following properties of $J:$
\begin{enumerate}
\item{the Zygmund classes $Z$ and $Z_0$ are preserved by
$J,$}
\item{$J$ is anti-involutory in the sense that $J^2=-I,$ and}
\item{$J(\sin kx) = \cos kx {\rm \ and \ }
J(\cos kx) = - \sin kx.$}
\end{enumerate}
Proofs of all of these statements are greatly simplified by
considering different possible extensions of $V$ to the complex
plane, as we shall see in chapter 2.

The anti-involutory property of $J$ yields an almost complex
structure on ${\cal Z}$ the tangent space to universal
Teichm\"uller space. In fact, whenever an anti-involutory
automorphism $J$ of a
 vector space $X$ over ${\mathbb R}$ is given, $X$ becomes a
vector space over ${\mathbb C}$ by defining, for every $v$ in $X,$
$$(a+ib)v=av
+ b J(v).$$ The reader should check that $\left(
(a_1+ib_1)(a_2+ib_2)
\right)v = (a_1 + i b_1)\left((a_2+ib_2)v \right).$

Instead of using the exponential map $z=e^{ix} \mapsto x$ as a
real-valued chart for the unit circle, one can use the
stereographic map $z \mapsto u$ from the circle to
$\overline{\mathbb R}$ where
\begin{equation}
\label{stereographic}
u = U(z)=
\frac{z+i}{iz+1}.
\end{equation}
  This map sends the four points
 $1,i,-1,-i$ on the unit circle to the four points $1,\infty,-1,0,$
 respectively, on the real axis.
The real-valued vector field $\tilde{V}(z)\frac{\partial}{\partial
z}$ on the circle $\{z:|z|=1\}$ is related to the vector field
$\hat{V}(u)$ defined for $u$ on the real axis by
$$\tilde{V}(z) =
\hat{V}(u) \left(\frac{-2}{(u+i)^2}\right).$$
Since we assume $\tilde{V}$ is continuous, and in particular
bounded on the circle, that implies at most
quadratic growth of $\hat{V}$ near $\infty.$
That is $$\hat{V}(u)=O(|u|^2)$$
as $|u| \rightarrow \infty.$

In the special case when $\tilde{V}(z)\frac{\partial }{\partial z}
=(c_0+c_1z+c_2 z^2) \frac{\partial}{\partial z}$ is a real-valued
M\"obius vector field, then $c_1$ is pure imaginary,
$c_2=-\overline{c_0},$
$$\hat{V}(u)=\left(c_0+c_1\left(\frac{i(u-i)}{(u+i)}\right)
+c_2\left(\frac{i(u-i)}{(u+i)}\right)^2\right)\left(\frac{-2}{(u+
i)^2}
\right)$$
is a quadratic polynomial in $u,$ and
$$V(x)=\frac{(c_0 +c_1e^{ix}-\overline{c_0}e^{2ix})}{ie^{ix}}=
\frac{1}{2} a_0 +a_1 \cos x + b_1 \sin x.$$
Here, $a_0, a_1,$ and $b_1$ are real and $a_0=2 {\rm \ Im \ } c_1,
a_1= 2 {\rm \ Im \ } c_2, b_1 = 2 {\rm \ Re \ } c_2.$

It will turn out that the quadratic polynomials are preserved by
the Hilbert transform and so $J$ is well-defined on the quotient
space
\begin{equation}
\label{Zyg}
{\cal Z}= \{\tilde{V} \in Z\}/\{ quadratic \
polynomials\}.
\end{equation}

In section 2.6, when we use complex methods to deal with
Hilbert transforms, we will find it useful to
map the interior of the circle to the upper half-plane by
the stereographic map $u=U(z)$ given in (\ref{stereographic}) and
then
compute the Hilbert transform
in the upper half-plane.  When this is done, it must be remembered
that the function $\hat{V}$ is permitted to have at most quadratic
growth near infinity.

\end{subsection}

\begin{subsection}{Scales and Trigonometric Approximation}
If a vector field $\tilde{V}(z)\frac{\partial}{\partial z},$
on the unit circle
is given by a finite sum of the form
\begin{equation}
\label{trig2}
\tilde{V}_n(z) = \sum_{k=-n}^{n} c_k z^k,
\end{equation}
and is real-valued,  then
$c_{n+2}= - \overline{c_{-n}}.$
The corresponding function
$V_n(x)=\tilde{V}_n(e^{ix})/ie^{ix},$
is the trigonometric polynomial
$$a_0/2 + \sum_{k=1}^{n}\left(a_k \cos kx + b_k \sin kx \right)$$
of degree $n,$ where
$a_k= 2 {\rm \
Im \ } c_{k+1}$ and $b_k= 2 {\rm \ Re \ }c_{k+1}.$
We may think of the trigonometric polynomial $V_n$ with
$$||V_n||_{\infty}=M$$
as a typical
vector field having a definite oscillation down to intervals whose
length is as small as $\frac{1}{M n}.$  That is, if $V(x)=1$
and  $0<t<\frac{1}{M n},$ then $V(x+t)>0.$
  This is because of the mean value theorem and the
following lemma due to Bernstein \cite{Bernstein}.
\begin{lemma}
\label{Bernstein}
  If $V_n(x)$ is a trigonometric polynomial of degree
$n,$ then
$$||\frac{d}{dx} V_n(x)||_{\infty} \leq n ||V_n(x)||_{\infty}.$$
\end{lemma}
\begin{proof}We follow the proof given in \cite{Lorenz}, page 39.
To begin, assume there is a trigonometric polynomial $V_n$ with
$||V_n'||=n L$ and $L>||V_n||.$  Thus at some point $x_0,$
$|V_n'(x_0)|=nL,$ and we can assume that $V_n'(x_0)=nL.$  Since
$V_n'$ is a maximum at $x_0,$ $V''(x_0)=0.$

Consider the trigonometric
polynomial
$$T_n(x) = L \sin n(x-x_0) - V_n(x)$$
of degree $n.$
In the interval $[x_0,x_0+2 \pi)$
there are $2n$ points
where $\sin n(x-x_0)$ takes the values $\pm 1,$  and between any
two of
these points the polynomial
$T_n$ takes values of opposite sign. Hence
$T_n$ has $2n$ different zeros in this interval, and
so
$$T_n'(x) = n L \cos n(x-x_0) - V_n'(x)$$
also has $2n$ different zeros.  One of these zeros is $x_0,$ since
$$T_n'(x_0)= n L - V_n'(x_0).$$
Also,
$$T_n''(x)=-n^2 L \sin n (x-x_0) - V_n''(x)$$
vanishes at $x=x_0.$  Moreover, $T_n''$ has $2n$ zeroes between the
zeros of $T_n'.$  Thus $T_n''$ has at least $2n+1$ zeros in this
interval, and since it is a trigonometric polynomial of degree $n$
it
must be identically zero.
Thus $T_n'$ is constant, but since $T_n'(x_0)=0,$ this implies
$T_n$ is
constant.  But this contradicts the statement that $T_n$ changes
sign
and we conclude that the original assumption could not be correct,
that
is, we have $L \leq ||V_n||_{\infty},$ which means
$||V_n'||_{\infty}
\leq n ||V_n||_{\infty}.$
\end{proof}

Assume $V$ and $W$ are continuous functions of period $2 \pi.$
An inequality of the form $||V(x)-W(x)||_{\infty}<1/2^n$ for large
$n$ implies that
the graph of $V(x)$ closely resembles the graph of $W(x).$
Now
consider $M_{k,I} V (x) = \frac{1}{2^k} V(2^k x),$
where $x$ lies in some interval $I$ of length
$2\pi/2^k.$
Then $M_{k,I}$ is a magnification operator of degree $k,$
magnifying the graph of $V$ over the interval $I$
by the same factor in both
the domain and range.

In fractal geometry, one considers graphs that have roughly
the same shape no matter how much
they are magnified.  Thus, suppose that we
go to some fine scale $M_{k,I}V.$  Then the picture of the graph
seen
at this scale should roughly resemble the picture of the graph of
$M_{n,J}$ if $n$ is any number larger than $k$ and $J$ is any
interval of size $2 \pi/2^n.$

A general trigonometric polynomial does not possess this property.
Suppose a polynomial $V_{2^n}$ of
degree $2^n$ is magnified by the operator $M_{k,I}$
of degree $k.$ If $k$ is larger
than $n,$ then because of Bernstein's inequality,
one does not see any oscillation in the graph of
$M_{k,I}V_n.$
This observation motivates the following theorem due to Zygmund
and Jackson, \cite{Zygmund}, \cite{Jackson}, \cite{Lorenz}.

\begin{thm}  Suppose $V(x)$ is a continuous, periodic function
defined on the real axis.  Then $V$ is in the Zygmund
class $Z$ defined by
$$\left| \frac{V(x+t)+V(x-t)-2V(x)}{t} \right| \leq C$$
if, and only if, there exists a constant $C'$ such that for every
positive integer $n$ there is a trigonometric polynomial $V_n$ of
degree at most $n,$ such that
$$||V-V_n||_{\infty} \leq \frac{C'}{n}.$$
Moreover, the number $C'$ can be estimated purely in terms of $C$
and vice versa.
\end{thm}
\begin{proof}
We begin by proving that if such trigonometric approximations are
possible for every $n,$ then $V$ is in the Zygmund class.  For each
integer of the form $2^k,$ let $V_{2^k}$ be a trigonometric
polynomial of degree $2^k$ such that
\begin{equation}
\label{geometric}
||V-V_{2^k}||_{\infty} \leq \frac{C}{2^k}.
\end{equation}

Now select $n$ so that $\frac{1}{2^{n+1}} \leq t \leq
\frac{1}{2^{n}},$
and write $V$ in the form $V=W_1+W_2$ where
$W_1=V-V_{2^n}$ and $W_2= V_{2^n}.$ In general, define the
difference
operator $\Delta_t$ by
$$\Delta_t G(x)=G(x+t)-G(x).$$
Then
$$\Delta^2_tG(x)=\Delta_t(\Delta_t(G))(x)=G(x+2t)-2G(x+t)+G(x).$$
From the hypothesis,
\begin{equation}
\label{doubleuone}
|\Delta^2_t W_1 (x)| \leq
\frac{8C}{2^{n+1}} \leq 8 C t.
\end{equation}

Putting $V_0=0,$ we may rewrite $W_2$
as a sum over scales:
\begin{equation}
\label{scalesum}
W_2(x) = V_1 - V_0 + \sum_{k=1}^n
\left(V_{2^k} - V_{2^{k-1}}\right).
\end{equation}
Each term
$V_{2^k}-V_{2^{k-1}}$ has norm bounded by
$$|V_{2^k}-V| +
|V-V_{2{k-1}}|
\leq \frac{3C}{2^k},$$ and is
a trigonometric polynomial of degree less than or equal to $2^k.$
So by Lemma \ref{Bernstein}, the second derivative of
$V_{2^k}-V_{2^{k-1}}$ is at most $3C 2^k.$ Thus, by
the second mean value theorem
$$ |\Delta^2_t \left(V_{2^k}-V_{2^{k-1}}\right)| \leq 3C 2^k t^2.$$
 By using equation (\ref{scalesum}), we obtain
 \begin{equation}
 \label{doubleutwo}
 |\Delta_t^2 W_2 | \leq
 \sum_{k=1}^{n} 3 C 2^k t^2 \leq \sum_{k=1}^{n}
 \frac{3 C 2^k}{2^{2n}} = \frac{3C 2^{n+1}}{2^{2n}}= \frac{6C}{2^n}
\leq
 12 C t.
 \end{equation}
 Putting inequalities (\ref{doubleuone}) and (\ref{doubleutwo})
 together, we obtain $|\Delta^2_t V| \leq 20 C t$ and this
 proves the first half of the theorem.

To prove the other half, for every $n$ we must construct a
trigonometric polynomial $V_n$ of degree $n$ that approximates $V$
in the sup-norm to within $C/n.$
Let $K_n$ be the Jackson kernel defined by
$K_n=\sigma_{2n-1} - 2
\sigma_n,$ where
$$\sigma_n=\frac{s_0+s_1+ \cdots +s_{n-1}}{n}
=\frac{1}{2\pi n} \left( \frac{\sin \frac{nt}{2}}{\sin \frac{t}{2}}
\right)^2$$
is the Fej\'er kernel and
$$s_n=\frac{ \sin (2n+1)\frac{t}{2}}{2 \sin \frac{t}{2}}.$$

By convolution of $V$ with the
Jackson kernel $K_n,$ one gets a trigonometric
polynomial of degree $2n-1$ that approximates
$V$ to within $O(1/n)$ in the sup norm.  For details of
this proof we refer to \cite{Lorenz}, pages 55-56, or to
\cite{Zygmund}.
\end{proof}

\begin{thm}
The Zygmund spaces $Z$ and $Z_0$ are invariant under the Hilbert
transform $J.$
\end{thm}
\begin{proof}
Since a much easier
proof of the same result is given in section 2.6, here we only
outline the argument given by Zygmund in \cite{Zygmund}.
Begin by using
a result of Favard \cite{Favard}:
if $|g'| \leq M$ then $|J g- J \sigma_n(g)| \leq A/n.$
Then employ the Zygmund-Jackson theorem.
Let $V$ be in $Z$ and $V_n$ be a trigonometric
polynomial of degree $n$ with $|V-V_n| \leq \frac{C}{n}.$
Let $G'=V$ and $T_n' = V_n.$  Then $|(G-T_n)'| \leq \frac{C}{n}$
and therefore  $|J (G-T_n) - J \sigma_n(G-T_n)| \leq
\frac{AC}{n^2}.$ Thus $JG$ is approximable in the sup norm to
within
$\frac{AC}{n^2}$ by a trigonometric polynomial of degree $n.$
This implies that $JV$ is approximable to within $\frac{AC}{n}$ by
a
trigonometric polynomial of degree $n,$ therefore $JV$ is in the
Zygmund class.

The proof for the class $Z_0$ is similar.
\end{proof}

\end{subsection}

\begin{subsection}{Automorphisms of \te\ Space}
Given any quasisymmetric homeomorphism  $f$ of $S^1,$ the map
$\rho_f([h])=[h \circ f^{-1}]$ is a bicontinuous self-map of
$T=QS {\rm \ mod \ }PSL(2,{\mathbb R}).$  Moreover, $\rho_f$
preserves the almost complex structure.  We call
biholomorphic automorphisms
of $T$ of this form {\em geometric} automorphisms.

An {\em almost complex
structure} on a real Banach manifold $M$ is a smoothly varying
family of automorphisms
$J_x, x \in M,$ of each fiber of the tangent bundle
such that $J_x^2 = - I.$ A diffeomorphism $F$ of $M$ is
{\em almost complex} if $J_{F(x)}\left(F'_x(v)\right) =
F'_x(J_x(v)).$

\begin{thm}
\label{autothm}
Any almost complex  automorphism  $F$ of $T$
is geometric.  That is, given a diffeomorphism $F$
of $$QS {\rm \ mod \ } PSL(2,{\mathbb R})$$
whose derivative commutes with the almost complex structure,
there exists a
quasisymmetric map $f$ such that $F=\rho_f.$
\end{thm}
An outline of the proof of this theorem is given in the last
section
of this paper.
\end{subsection}
\end{section}

\begin{section}{Complex Analysis}

\begin{subsection}{Quasiconformal Extensions}
Roughly speaking, a homeomorphism
of $\mathbb{R}^n$ is quasiconformal
if it distorts standard shapes by a bounded amount, see \cite{LV},
\cite{Rickman},
\cite{GL}. When $n \geq 2,$ it turns out that quasiconformal maps
are
differentiable almost everywhere and the distortion of
shape can be measured infinitesimally.  An observation of central
importance
for the deformation theory of one-dimensional dynamical systems is
that
this statement is not true when $n=1.$  That is, quasisymmetric
maps may not be differentiable anywhere.

In any case, the measurement of quasiconformal distortion at a
point $z$ for
a mapping $f$ when $n=2$
is by means of a quantity called the {\em local dilatation}
$K_z(f):$
\begin{equation}
\label{Dil}
K_z(f) = \frac{|f_z|+|f_{\overline{z}}|}{
|f_z|-|f_{\overline{z}}|}.
\end{equation}

Any quasisymmetric homeomorphism $h$ of the real axis extends to
a
quasiconformal self-mapping of the upper half-plane.  This
pivotal result was first proved by Ahlfors and Beurling \cite{BA}.
The
formula given in \cite{A} for such an extension
of $h$ is
$H_1(z)=F(z) + i G_1(z),$ where
\begin{equation}
\label{BA}
\begin{array}{l}
F(x+iy)= \frac{1}{2y} \int_{x-y}^{x+y} h(t) dt \\
\ \\
G_1(x+iy)= \frac{1}{2y} \left\{\int_x^{x+y} h(t)dt -
\int_{x-y}^y h(t) dt \right\}.
\end{array}
\end{equation}
This formula does not
extend the identity by the identity.  In particular, for $h(x)=x,$
the extension $H_1(z)=x + \frac{1}{2} i y.$  It is therefore
convenient to
multiply the expression for $G_1$ by a factor two.  That is, we put
\begin{equation}
\label{BA1}
G(x+iy)= \frac{1}{y} \left\{\int_x^{x+y} h(t)dt -
\int_{x-y}^y h(t) dt \right\}.
\end{equation}
Although $H=F+iG$ differs from $H_1,$ it still yields
a quasiconformal extension $ex(h)$
of $h,$ but with the additional property that the identity is
extended by the identity.
 It is useful
to view $H$ as an extension to the whole plane by stipulating
that $H(\overline{z})= \overline{H}(z).$  The reader should check
that this extension process is
natural for real affine transformations in the sense that if
$A(z)=c_1z+c_2$ and $B(z)=c_3z + c_4$ where $c_1, \ldots, c_4$ are
real, then
$$ex(A \circ h \circ B)= A \circ ex(h) \circ B.$$
Hence, if we assume
\begin{equation}
\label{lift}
h(x)+1=h(x+1),
\end{equation}
 then $H(z+1)=H(z)+1.$

It is important to note that a lift of a self-homeomorphism of the
circle by the
universal covering $x \mapsto e^{2 \pi i x}$ yields a
homeomorphism $h$ satisfying (\ref{lift}), and conversely,
if a homeomorphism of the real-axis satisfies (\ref{lift}) then it
projects to a homeomorphism of the circle.  Moreover,
the covering $x \mapsto e^{2 \pi i x}$ extends to the covering
$z \mapsto e^{2 \pi i z}$ of the punctured unit disc
${\mathbb D}^*={\mathbb D}-\{0\}$ by
the upper half-plane.  The extension of $h$ to the disc punctured
at $0$
is a quasiconformal map preserving $0,$ and therefore if we
stipulate
that the extension
preserves $0,$ it becomes a quasiconformal extension to the entire
disc.
This method of extension also has the following {\em asymptotic
property}:-

Assume that $h(0)=0, h(x+1)=h(x)+1,$  $h$ is
quasisymmetric and
$$\frac{1}{1+\epsilon} \leq \frac{h(x+t)-h(x)}{h(x)-h(x-t)} \leq
1+\epsilon$$
for $|t|<\epsilon,$ with $\delta$ sufficiently small.
Then if $| {\rm Im \ } z |< \delta$,
the dilatation $K_z$ of $ex(h)$ at $z$
satisfies $K_z <1+\epsilon',$ where
$\epsilon'$ converges to zero as $\epsilon$ converges to zero.

\end{subsection}

\begin{subsection}{\tes\ Metric}
The \te\ distance between two points $[h_1]$ and $[h_2]$ in
$QS {\rm \ mod \ } PSL(2,\mathbb{R})$ is defined to be
\begin{equation}
\label{tdistance}
d([h_1],[h_2]) = \frac{1}{2} \log K_0(h_2 \circ (h_1)^{-1}),
\end{equation}
where
$$K_0(h) =
 \inf \{ K(\tilde{h}): {\rm \ where \ } \tilde{h} {\rm \ is \ any
\ quasiconformal \ extension \ of \ } h \}.$$

As a consequence of basic properties of quasiconformal mappings,
$d([id],[h])$ is
always realized by an extremal mapping $\tilde{h}$ which is an
extension
of $h.$  To see this one can assume that $h$ fixes three points on
the real axis, say $0,1$ and $\infty,$ and then select extensions
$\tilde{h}_n$ of $h$ such that $\frac{1}{2} \log K(\tilde{h}_n) <
\frac{1}{2} \log K_0([h])+\frac{1}{n}.$ Since the mappings
$\tilde{h}_n$
are normalized and have uniformly bounded dilatation, they are
equicontinuous.  Therefore there is a subsequence of $\tilde{h}_n$
that
converges uniformly in the spherical metric to some self-mapping
of the
upper half-plane $\tilde{h}_0.$ Since each $\tilde{h}_n$ coincides
with
$h$
at every point of the real axis, so does $\tilde{h}_0.$  Moreover,
the
maximal
dilatation of $\tilde{h}_0$ must be less than $K([h])+\frac{1}{n}$
for
every positive integer $n.$  On the other hand $K(\tilde{h}_0)$
cannot
be less
than $K_0([h])$ because by definition $K_0([h])$ is the infimum of
the
dilatations
of all  possible extensions of $h.$ We conclude that
every mapping $h$ of the real axis has an {\em extremal}
quasiconformal extension $\tilde{h}_0$ to the upper half-plane,
that is,
an extension for which
$$K_0([h])=K(\tilde{h}_0).$$

It will turn out that certain quasisymmetric mappings $h$ have many
extremal extensions and thus we do not expect to find a general
formula
that yields an extremal extension. In particular, the
Beurling-Ahlfors
extension formula
given in the previous section will almost never yield an extremal
extension.

In formula (\ref{tdistance}) the \te\ distance is seen as the
solution
to an infimum problem. It turns out that it is also the solution
to a
supremum problem.  Consider the vector space  ${\cal R}({\mathbb
H})$
of integrable holomorphic quadratic
differentials $\varphi(z)dz^2$ in the upper half-plane ${\mathbb
H}$
with only
a finite number of poles on the real axis and for  which
$\varphi(z)dz^2$
is real-valued on the real axis. Any element
$\varphi(z)dz^2$ of ${\cal R}({\mathbb H})$ has the form
\begin{equation} \label{rationalqds}
\varphi(z)dz^2 = \frac{p(z)dz^2}{(z-x_1)\cdots(z-x_n)},
\end{equation}
where $x_1,\ldots,x_n$ are distinct points on the real axis and
$p(z)$ is a polynomial of degree less than or equal to $n-3$ with
real
coefficients.  Such a quadratic differential determines a
decomposition
of ${\mathbb H}$ into a finite number of strips, $S_1, \ldots,
S_k,$
where $k \leq n-3.$ The
interior of each strip is swept out by horizontal trajectories of
this
quadratic differential, that is, parameterized curves $\alpha(t)$
along
which $\varphi(\alpha(t))\alpha'(t)^2 >0.$  A choice of coordinate
$$\zeta = \pm \int \sqrt{\varphi(z)} dz + (const)$$
can be made so that $z \mapsto \zeta$ maps
the $j$-th strip to a rectangle $R_j$ and takes the horizontal
trajectories
$\alpha(t)$ to horizontal line segments that join the left side of
the
rectangle to its right side.
Let $a_j$ and $b_j$ be the width and height of the rectangle $R_j$
measured in the parameter $\zeta.$
Then
$$||\varphi||=\int \! \int_{\mathbb H} |\varphi(z)| dx dy =
\sum_{j=1}^k
a_j b_j.$$ That is $||\varphi||$ is equal to the sum of the areas
of
these rectangles.
Moreover, if $\beta$ is any arc in ${\mathbb H}$ with endpoints on
${\partial} {\mathbb H}$ transversal to the horizontal trajectories
of
$\varphi,$ then we can assign to it a {\em height}
$ht_\varphi(\beta)$
given by
$$ht_{\varphi}(\beta) = \int_{\beta} {\rm \ Im}(\sqrt{\varphi(z)}
dz),$$
equal to the sum of the heights of the rectangles corresponding to
the strips $S_j$ crossed by $\beta.$

Let $I_j$ be the intervals on $\partial {\mathbb H}$ whose
endpoints are successive pairs from
the (ordered) sequence of points $x_1, \ldots , x_k$ and assume
that the endpoints of the arc
$\beta\subset {\mathbb H}$
lie on the intervals $I_{j_1}$ and $I_{j_2}.$  If $\beta'$ is
another
arc transverse to the horizontal trajectories of $\varphi(z)dz^2$
with
endpoints on the same intervals $I_{j_1}$ and $I_{j_2},$ then
$ht_{\varphi}(\beta)=ht_{\varphi}(\beta'),$ and so the height
function
$ht_{\varphi}$ is a nonnegative function defined on all possible
pairs
of intervals $I_{j_1}$ and $I_{j_2}$ taken from the set $I_1,
\ldots ,
I_k.$

A sense-preserving selfmap $h$ of $\partial {\mathbb H}$ takes
the
points $x_j$ to points $x_j'=h(x_j)$ and thus determines a new
height
function defined on pairs of intervals $I_{j_1}'=h(I_{j_1})$ and
$I_{j_2}'=h(I_{j_2}).$
From the theorem of Hubbard and Masur (in \cite{HM};see also
\cite{G1})
there is a unique quadratic differential of the form
$$\psi(z)dz^2 = \frac{q(z)dz^2}{(z-x_1')\ldots(z-x_n')},$$
such that q(z) has real coefficients and degree $\leq n-3$
and the heights of $\psi$ for the interval pairs $I_{j_1}'$ and
$I_{j_2}'$ are
equal to the heights of $\varphi$ for the corresponding interval
pairs
$I_{j_1}$ and $I_{j_2}.$
Moreover, $\psi$ is unique among all continuous integrable
quadratic
differentials on $\overline{\mathbb C} \setminus
\{x_1',\ldots,x_k'\},$
real-valued on the real axis, with heights between pairs of
intervals $I_{j_1}'$ and $I_{j_2}'$ greater than or equal to the
corresponding heights of $\varphi$ on intervals $I_{j_1}$ and
$I_{j_2}$
and with $||\psi||$ as small as possible.

If $h$ has a quasiconformal extension $\tilde{h}$ with dilatation
$K_0,$
then $||\psi|| \leq K_0 ||\varphi||,$ and one obtains the following
expression for $K_0:$
\begin{equation}
\label{tdsup}
K_0 = \sup \frac{||\psi||}{||\varphi||},
\end{equation}
where the supremum is taken over all non-zero quadratic
differentials
$\varphi$ of the
form (\ref{rationalqds}) and $\psi$ is the quadratic differential
with
simple poles at the points $x_j'=h(x_j), 1 \leq j \leq n,$ and with
the
same corresponding
heights with respect to these points that $\varphi$ has with
respect to
the points $x_j, 1 \leq j \leq n.$
 \end{subsection}

\begin{subsection}{Quadratic Differentials}
Let ${\cal A}={\cal A}(\Omega)$ be the Banach space of integrable
functions $\varphi(z)$ holomorphic in ${\Omega}$ where $\Omega=
{\mathbb H}$ or $\Omega=\Delta=\{z:|z|<1\}$
with norm
$$||\varphi|| = \iint_{\Omega} |\varphi(z)|dxdy < \infty.$$
In this section we introduce a pairing between ${\cal A}$ and
${\cal Z}$ and show that ${\cal A}^* \cong {\cal Z}$
and that $({\cal Z}_0)^{*} \cong {\cal A}.$  By ${\cal Z}$ we mean
the
vector fields $V$ defined on $\partial \Omega$ such that
$V(z)\frac{\partial}{\partial z}$ is real-valued on $\partial
\Omega,$
and such that $||V||_{\cal Z} < \infty.$
Since there is a M\"obius transformation transforming ${\mathbb H}$
onto
$\Delta$ and since the statements we prove will be invariant under
pull-back by M\"obius transformations, we can work interchangeably
with
either ${\mathbb H}$ or with $\Delta.$

Our first step is to prove a special
case of Bers' approximation theorem \cite{Bers1}, \cite{A1}.
Let ${\cal R}({\mathbb H})$
be the space of finite linear combinations of the
form
\begin{equation}
\label{lincomb}
\lambda_1 \varphi_{x_1}(z) + \cdots + \lambda_n \varphi_{x_n}(z),
\end{equation}
where $x_1, \ldots, x_n$ and $\lambda_1, \ldots, \lambda_n$ are
real
numbers and
$$\varphi_x(z) = \frac{x(x-1)}{z(z-1)(z-x)}.$$

 \begin{thm}
\label{BAT}
${\cal R}$ is dense in
${\cal A}.$ \end{thm}
\begin{proof}
A similar and much deeper result is true if ${\mathbb H}$ is
replaced by any plane domain, \cite{A1}, \cite{Bers1}.

Let $L$ be any linear functional on the Banach space ${\cal A}$
that
annihilates ${\cal R}.$  To show that ${\cal R}$ is dense in ${\cal
A}$ it is sufficient to show that $L$ annihilates ${\cal A}.$ By
the
Hahn-Banach and Riesz representation
theorems, there exists a bounded measurable function $\mu$ defined
in
${\mathbb H}$ so that
$$L(\varphi) = {\rm real \ part \ of }\left( \iint_{\mathbb H}
\mu(z) \varphi(z) dx dy \right).$$
If we extend $\varphi(z)$ and $\mu(z)$ to the lower half-plane by
the
rules $\overline{\varphi(z)}=\varphi(\overline{z})$ and
$\overline{\mu(z)}=\mu(\overline{z}),$
then we can write the formula for $L$ as
$$L(\varphi) = \frac{1}{2} \iint_{\mathbb C} \mu(z) \varphi(z) dx
dy.$$

The assumption that $L$ annihilates ${\cal R}$ implies that
\begin{equation}
\label{vf1}
V(z) = - \frac{z(z-1)}{\pi} \iint_{\mathbb C}
\frac{\mu(\zeta)}{\zeta(\zeta-1)(\zeta-z)} d\xi d\eta =0
\end{equation}
whenever $z$ is a real number.
One shows that $V$ has the following properties:
\begin{enumerate}
\item{$\overline{\partial} V
= \mu$ in the sense of distributions,}
\item{ $|V(z)|=O(|z| \log |z|)$ as $z \rightarrow \infty,$ and}
\item{$V(z)$ has an $|\epsilon \log \epsilon|$-modulus of
continuity,
that is to say, given $R>0,$ there exists a $C$ such that for every
$z_1$ and $z_2$ with $|z_1|$ and $|z_2|<R$ and with $|z_1-z_2| <
1/2,$
$$|V(z_1)-V(z_2)| \leq C |z_1-z_2| \log (1/|z_1-z_2|).$$}
\end{enumerate}

Let
$D_{\epsilon,R}$ be a semi-disc in the
upper half-plane with
diameter of
length $2R$ along the line $y=\epsilon$ and with
midpoint on the $y$-axis. The curved part of the boundary of
$D_{\epsilon}$ is parameterized by the curve $z=i\epsilon +
Re^{i \theta}, 0 \leq \theta \leq \pi.$
Assume further that
$\varphi$ is continuous on the real axis and $\varphi(z) =
O(|z|^{-4})$ as $z \rightarrow \infty.$ Then
the subspace of $A({\mathbb H})$ comprising those $\varphi$ with
these properties is dense in $A({\mathbb H}).$
Since $\varphi$ is integrable and $\mu$ is
bounded, $$\iint_{\mathbb H} \mu \varphi = \lim
\iint_{D_{\epsilon, R}} \mu d \xi d \eta,$$
where the limit is taken both as
$\epsilon \rightarrow 0$ and as
$R \rightarrow \infty.$
On the other hand, from Green's formula,
$$\iint_{D_{\epsilon,R}} \mu d \xi d \eta= \int_{\partial
D_{\epsilon,R}} V(\zeta) \varphi(\zeta) d \zeta.$$
Because $V(z)$ is identically zero when $z \in
{\mathbb R},$ if we first take the limit in this line integral  as
$\epsilon \rightarrow 0$ we obtain
$$
\int_0^{\pi} V(Re^{i\theta})
\varphi(Re^{i\theta})Rie^{i\theta} d
\theta = 0.$$
Because of the vanishing condition on $\varphi,$
$$\int_0^{\pi}
V(Re^{i\theta})
\varphi(Re^{i\theta})Rie^{i\theta} d
\theta $$
is dominated by a constant times $(\log R)/R$ and thus vanishes as
$R
\rightarrow \infty.$
 \end{proof}

We are now ready to introduce the pairing between an element $V$
in
${\cal Z}$ and $\varphi$ in ${\cal A}.$
Given $V$ in ${\cal Z}$ we select any extension $\tilde{V}$ of $V$
to
the upper half-plane with the properties that $\partial \tilde{V}
=
\mu$ is essentially bounded and that $|V(z)|=O(|z|^2).$  Then we
define
\begin{equation}
\label{pairing}
(V,\varphi) = {\rm Re} \left( \iint_{\mathbb H} \mu
\varphi \right).
\end{equation}

We must show first that any $V$ in ${\cal Z}$ has such an
extension and second that if a different extension is taken, the
integration (\ref{pairing}) yields the same result.
$\tilde{V}$ can be defined by the Beurling-Ahlfors' formula
(\ref{BA}) applied to the vector field $V$:
$$ {\rm Re }(\tilde{V}(x+iy))
= \frac{1}{2y}\int_{x-y}^{x+y}V(t) dt {\rm \ \  and}$$
$$ {\rm Im }(\tilde{V}(x+iy)) =
\frac{1}{y} \left( \int_{x}^{x+y}V(t) dt -
\int_{x-y}^{x}V(t) dt \right).$$
We leave it to the reader to verify that $\tilde{V}$ has the
appropriate
growth rate and that the Zygmund condition implies
$$\overline{\partial} \tilde{V} =
\frac{1}{2}(\frac{\partial}{\partial x} + i
\frac{\partial}{\partial y})
\tilde{V}$$
is bounded.
To show that the right hand side of (\ref{pairing}) depends only
on the
values of $V$ on the real axis we first note that from Theorem
\ref{BAT}
it suffices to show the right hand side of (\ref{pairing}) depends
only
on the values of $V$ on the real axis when $\varphi$ has the
special
form (\ref{lincomb}). In that case, if we assume $V$ vanishes at
$0$ and
$1$ and has growth rate $o(|z|^2),$ then by Green's formula,
\begin{equation}
\label{residue}
(V, \varphi) = \frac{\pi}{2} \sum_j \lambda_j V(x_j).
\end{equation}

\begin{thm} The pairing (\ref{pairing}) induces an isomorphism
between
${\cal Z}$ and ${\cal A}^*.$
\end{thm}
\begin{proof}
From the preceding discussion we have seen how an element $L$ of
${\cal A}^*$ determines by the correspondence an element $V$ in
${\cal
Z}.$ Conversely, because of
the residue formula (\ref{residue}), the extension formula and
because
${\cal R}$ is dense in ${\cal A}$, any element $V$ of ${\cal Z}$
determines by this correspondence an element of
${\cal A}^{*}.$
Note that $||V||_{cr}=0$ is equivalent to the condition that
$V(z)=a_0+a_1z+a_2z^2.$ Since by definition elements
$$\varphi(z) = \sum_{j=1}^n \frac{\lambda_j}{z-x_j}$$
of ${\cal R}$ satisfy $|\varphi(x)|=O(|z|^{-4})$ as $z \rightarrow
\infty,$ they also satisfy $\sum
\lambda_j =0, \sum_j \lambda_j x_j =0$ and $\sum_j \lambda_j
x_j^2=0.$
One therefore sees that
any quadratic polynomial vector field $V(x)
\frac{\partial}{\partial x}$
annihilates all elements of ${\cal R},$
and so also annihilates ${\cal A}$ since
${\cal R}$ is dense in ${\cal A}.$
\end{proof}

Because ${\cal A}^*$ is isomorphic to ${\cal Z},$ the norm on
${\cal Z}$
dual to the norm on ${\cal A}$ is equivalent to $|| \ ||_{cr}.$
We call this norm the {\em infinitesimal \te\ norm} and denote it
by
$|| \  ||_T.$  It is given by either of the following formulas:
$$||V||_T= \sup
\left\{ \left| \iint_{\mathbb H} \varphi \overline{\partial}
\tilde{V}
dx dy \right|: \varphi \in {\cal A} {\rm \ with \ }
||\varphi||=1 \right\}
$$
$$ = \inf \{||\overline{\partial} \tilde{V}||_{\infty}:
{\rm \ where \ } \tilde{V} {\rm \ is \ any \ extension \ of \ } V.
\}$$

\begin{definition}
We  say a sequence $\varphi_n$ in ${\cal A}$ is {\em degenerating}
if
there is a constant $C>1$ such that $C^{-1} \leq ||\varphi_n|| \leq
C$
and $\varphi_n(z) \rightarrow 0$ for every $z \in
{\mathbb H}.$
\end{definition}

We now wish to focus attention on the closed subspace ${\cal Z}_0$
of
${\cal Z}$ defined in section 1.6, the tangent vector fields to the
symmetric circle maps.
\begin{thm}
The following conditions on an element $V$ of ${\cal Z}$ are
equivalent:
\begin{enumerate}
\item{with respect to any smooth local coordinate $x$ on the
boundary of
$\Omega,$
$$ \left|\frac{V(x+t)-2V(x)+V(x-t)}{t}\right| \leq c(t),$$
where $c(t)$ approaches $0$ as $t \rightarrow 0,$}
\item{$V$ has a continuous extension $\tilde{V}$ for which
$\overline{\partial} \tilde{V} = \mu$ is vanishing in the sense
that
for every $\epsilon>0$ there exists a compact subset of ${\mathbb
H}$
such that if $z$ lies outside the compact set then $|\mu(z)| <
\epsilon,$}
\item{$V$ annihilates every degenerating sequence in ${\cal A},$}
 \end{enumerate}
\end{thm}
\begin{proof}  Given $V$ satisfying condition {\em 1}, the
Beurling-Ahlfors
formula
(\ref{BA}) yields a vector field $\tilde{V}$ with the property that
$\overline{\partial} \tilde{V}$ is vanishing.  Thus {\em 1} implies
{\em 2}.
It is easy to see that if $|\mu(z)|<\epsilon$ for $z$ outside a
sufficiently large compact set and if $\varphi_n$ is degenerating,
then
$$\lim_{n \rightarrow \infty} \iint \mu \varphi_n dx dy \rightarrow
\infty,$$
and so {\em 2} implies {\em 3}.
To see that {\em 3} implies {\em 1}, consider the following
sequence of
quadratic
differentials $\varphi_n,$ where $t_n \rightarrow 0$ and $x_n$ is
arbitrary:
 $$\varphi_n(z) = \frac{1}{t_n} \left\{ \frac{1}{z-(x_n-t_n)}
-\frac{2}{z-x_n}+ \frac{1}{z-(x_n+t_n)} \right\}$$
$$=
\frac{2t_n}{(z-(x_n-t_n))(z-x_n)(z-(x_n+t_n))}.$$
Note that $$\iint_{\mathbb H} |\varphi_n| = \iint_{\mathbb H}
\left|\frac{2t_n}{(z-t_n)z(z+t_n)}\right| dx dy =
\iint_{\mathbb H} \left|\frac{1}{(z-1)z(z+1)}\right| dx dy,$$
which is a positive constant not depending on $n$ and, for fixed
$z$ in
the upper half-plane, $\varphi_n(z) \rightarrow 0$ as $t_n
\rightarrow
0.$

By formula (\ref{residue})
$$(V,\varphi_n)= \frac{\pi}{2} \left\{
\frac{V(x_n-t_n)-2V(x_n)+V(x_n+t_n)}{t_n} \right\}$$
and we know that this quantity approaches zero as $t_n \rightarrow
0,$
no matter which sequence $\{x_n\}$ is
selected.  Thus {\em 3} implies {\em 1}.
\end{proof}

\begin{thm} The pairing $(V,\varphi)$ defined in (\ref{pairing})
induces
an isomorphism from ${\cal A}$ onto ${\cal Z}_0^*.$
\end{thm}
\begin{proof} We first observe that the pairing
defined in (\ref{pairing}) is non-degenerate
between ${\cal Z}_0$ and ${\cal A}.$  Since it is a non-degenerate
pairing between $A$ and ${\cal Z}$ and since ${\cal Z}_0 \subset
{\cal
Z},$  whenever $V \in {\cal Z}_0,$ $(V,\varphi)=0$ for all $\varphi
\in
{\cal A}$ implies $V=0.$  Moreover, for $|z_0|<1$ and
$\epsilon<1-|z_0|,$ by the mean value property
$$\varphi(z_0) = \frac{1}{\pi \epsilon^2} \iint_{|z-z_0|<\epsilon}
\varphi(z) dx dy= \iint \mu_0 \varphi,$$
where
$$
\mu_0(z) =\left\{
\begin{array}{ll}
\frac{1}{\pi \epsilon^2} & {\rm \ for \ } |z-z_0| \leq \epsilon,
\\
\ & \ \\
0 & {\rm \ for \ } z {\rm \  elsewhere}.
\end{array}
\right.
$$

Because the pairing is non-degenerate, the mapping from ${\cal A}$
to
${\cal Z}_0^*$ given by $\varphi \mapsto \{ V \mapsto (V,\varphi)
\}$
is well-defined and injective.  In order to show it is surjective
it
suffices to show
the unit ball of ${\cal A}$ is compact with
respect to the weak topology.  To this
end,  assume $\varphi_n$ is a sequence in ${\cal A}$ with
$||\varphi_n|| = 1$ and $L$ is a linear functional of the form
$$L(\varphi)= \iint_{\mathbb H} \mu \varphi,$$
where $|\mu(z)|<\epsilon$ for $z$ outside sufficiently large
compact
subsets of ${\mathbb H}.$  By normal convergence $\varphi_n$ has
a
subsequence that converges uniformly on compact subsets to some
$\varphi,$ which (in order to avoid cumbersome notation) we denote
also by $\varphi_n.$
Note that by Lebesgue dominated convergence
$$\lim_{n \rightarrow \infty} \iint
(|\varphi_n-\varphi|-|\varphi_n|)
= \lim_{n \rightarrow
\infty}(|\varphi_n-\varphi|-|\varphi_n|)=-||\varphi||,$$
and, since $||\varphi_n||=1,$ $||\varphi_n-\varphi|| \rightarrow
1-||\varphi||.$

We divide the argument into three cases: either
$||\varphi||=0$ or $0 < ||\varphi|| <1$ or $||\varphi||=1.$ In the
first case, $\varphi_n$ is degenerating
and $L(\varphi_{n}) \rightarrow 0.$
In the second case,
since $||\varphi||<1,$ if we put
$$\tilde{\varphi}_n = \frac{\varphi_n -
\varphi}{||\varphi_n-\varphi||},$$
then the denominator is bounded away from zero.  Thus
$\tilde{\varphi_n}$ is a degenerating sequence and
$L(\tilde{\varphi}_n)$ converges
to zero, which implies $L(\varphi_n)$ converges to
$L(\varphi).$ In the third case
$||\varphi_n - \varphi|| \rightarrow
\infty,$ and so
$L(\varphi_n)$ converges to $L(\varphi).$  Thus, in all cases
$L(\varphi_n)$ converges to $L(\varphi).$
 \end{proof}

\begin{definition}
We say a sequence $V_n$ in ${\cal Z}$ is {\em
vanishing} if there is a constant $C>1$ such that
$C^{-1} \leq ||V_n||_{T} \leq C$ and $V_n(x)$ approaches zero for
every
$x$ on the boundary.
\end{definition}

The following theorem enables one to deduce that an automorphism
of
${\cal Z}$ that is an isometry for the infinitesimal \te\ norm
necessarily preserves the closed subspace ${\cal Z}_0.$  This is
a key
step in the proof of the automorphism theorem, Theorem
\ref{autothm}
of section 1.8.

\begin{thm}
\label{thm11}
An element $V$ of ${\cal Z}$ with $||V||_T=1$ is in ${\cal
Z}_0$ if, and only if, for every vanishing sequence of elements
$W_n$ in
${\cal Z}$ with $||W_n||_T=1,$
$$\limsup_{n \rightarrow \infty} ||V+W_n||_T \leq ||V||_T.$$
\end{thm}
For the proof we refer to \cite{EG} and to \cite{GL}.

\end{subsection}

\begin{subsection}{Reich-Strebel Inequalities}
Let $\mu$ be the Beltrami
coefficient of a quasiconformal extension $\tilde{h}$ of a
quasisymmetric mapping $h,$ and let $K_0=K_0(h)$ be the smallest
possible dilatation of a quasiconformal extension of $h.$  For any
holomorphic quadratic differential $\varphi$ in ${\cal A}({\mathbb
H})$
with $$||\varphi||=\int \! \int_{\mathbb H}
|\varphi(z)| dx dy =1,$$ one has the following bounds on $K_0:$
\begin{equation}
\label{f1}
\frac{1}{K_0} \leq \int \! \int_{\mathbb H}
\frac{|1-\mu \frac{\varphi}{|\varphi|}|^2}{1-|\mu|^2}
|\varphi(z)|dx
dy,
\end{equation}
and
\begin{equation}
\label{f2}
K_0 \leq
\sup_{||\varphi||=1}
\int \! \int_{\mathbb H}
\frac{|1+\mu \frac{\varphi}{|\varphi|}|^2}{1-|\mu|^2}
|\varphi(z)|dx
dy.
\end{equation}
These inequalities were proved by Reich and Strebel \cite{RS1},
\cite{RS2}, who also observed that they yield the infinitesimal
form of
\tes\ metric:
\begin{equation}
\label{infform}
d_T([0],[t\mu]) = \frac{1}{2} \log K_0(t\mu) = t
\sup_{||\varphi||=1}
\left| \int \! \int_{|z|<1} \mu(z) \varphi(z) dx dy \right| + o(t),
\end{equation}
for $t>0.$
The inequality
$$d_T([0],[t\mu]) \geq t {\rm \ Re \ }
\left\{ \iint_{|z|<1} \mu \varphi \right\}
+ o(t),$$ for $||\varphi||=1,$ follows by replacing $\mu$ by $t\mu$
and
calculating the first variation in (\ref{f1}).
Similarly, the inequality
$$d_T([0],[t\mu]) \leq t \sup_{||\varphi||=1} {\rm \ Re \ }
\left\{ \iint_{|z|<1} \mu \varphi \right\} + o(t),$$
follows on replacing $\mu$ by $t\mu$ and calculating the first
variation in (\ref{f2}).

\end{subsection}

\begin{subsection}{Tangent Spaces Revisited}
Let $M$ denote the open unit ball in $L_{\infty}({\mathbb H})$ and
suppose $\mu_t$ is a smooth curve in $M$ such that each
$f^{\mu_t}$ is equal
to the identity on the boundary of ${\mathbb H}.$  Then from
inequality
(\ref{f1}) we see that for every holomorphic quadratic differential
$\varphi(\zeta)d\zeta^2$ on ${\mathbb H},$
$$1 \leq \iint_{\mathbb H}
\frac{|1-\mu_t \frac{\varphi}{|\varphi|}|^2}{1-|\mu_t|^2}
|\varphi| d \xi d \eta.$$
By putting
\begin{equation}
\label{orth}
||\mu_t -t \nu||_{\infty} = o(t)
\end{equation}
and computing the first variation in this inequality, one obtains
\begin{equation}
\label{orth1}
\iint_{\mathbb H} \varphi(\zeta) \nu(\zeta) d \xi d \eta = 0,
\end{equation}
for every such $\varphi.$

Conversely, suppose (\ref{orth1}) holds for every holomorphic
quadratic differential $\varphi$ on the upper half-plane.
Restricting $\varphi(\zeta) d \zeta^2 = \frac{d
\zeta^2}{(\zeta-z)^4}$
to the lower half-plane, we conclude
that
 \begin{equation}
\label{orth2}
\iint_{\mathbb H} \frac{\nu(\zeta)}{(\zeta-z)^4} d \xi d \eta =0,
\end{equation}
where $\zeta=\xi + i \eta$ is in the upper half-plane.

A key existence theorem (see \cite{G}, page 107) says that
(\ref{orth2})
implies there exists a curve $\mu_t$ such that
$$||\mu_t -t \nu||_{\infty} = o(t)$$
and $f^{\mu_t}(z)=z$ for every $z$ in the closure of the lower
half-plane.  Here, $f^{\mu_t}$ is the unique
quasiconformal
self-map of the whole plane that fixes $0,1$ and $\infty,$ and that
has
Beltrami coefficient equal to $\mu_t$ in the upper half-plane and
identically equal to $0$ in the lower half-plane.  In particular,
if we
let $M_0$ be the closed subspace of those $\mu$ in  $M$ for which
$f^{\mu}(x)=x$ for all $x \in {\mathbb R},$ then the tangent space
$N$
to $M_0$ consists of those $\nu$ for which
$$\iint_{\mathbb H} \nu(\zeta) \varphi(\zeta) d \xi d \eta =0$$
for every  quadratic differential $\varphi$ holomorphic in
${\mathbb
H}.$
Moreover, the tangent space to \te\ space $T$ is isomorphic to the
factor space
$L_{\infty}({\mathbb H})/N,$ see \cite{G} and \cite{GL}

For $\mu$ in $L_{\infty}(\mathbb H),$ define
$$\hat{\mu}= \left\{
\begin{array}{ll}
\mu(\zeta) & \mbox{ for $\zeta$ in \ } {\mathbb H} \\
\ & \ \\
\overline{\mu(\overline{\zeta})} & \mbox{ for $\zeta$ in \ }
{\mathbb
H}^*.
 \end{array}
\right.
$$
Let $\alpha:L^{\infty}({\mathbb H}) \rightarrow
Z$ be the map $\alpha: \mu \mapsto
V_{\mu}(x)$ where $x$ is in ${\mathbb R}$ and
 \begin{equation}
\label{vf}
V_{\mu}(z) = - \frac{1}{\pi} \iint_{\mathbb C}
\frac{\hat{\mu}(\zeta)}
{\zeta(\zeta-1)(\zeta-z)} d \xi d \eta.
\end{equation}
Also define
 $$\tilde{\mu}(\zeta) =
\left\{
\begin{array}{ll}
\mu(\zeta) & \mbox{ for $\zeta$ in \ } {\mathbb H}
\\ \ & \ \\
0 & \mbox{ for $\zeta$ in \ } {\mathbb H}^*,
\end{array}
\right.
$$
and let $\beta$ be the Bers map $\beta: \mu \mapsto \left(
W_{\mu}\right)'''(z),$ where
 \begin{equation}
\label{bcf}
W_{\mu}(z) = - \frac{1}{\pi} \iint_{\mathbb C}
\frac{\tilde{\mu}(\zeta)} {\zeta(\zeta-1)(\zeta-z)} d \xi d \eta.
\end{equation}

\begin{thm}
The maps $\alpha$ and $\beta$ defined above induce  isomorphisms
of
Banach spaces
from $L_{\infty}({\mathbb H})/N$ onto ${\cal Z}$
and from $L_{\infty}({\mathbb H})/N$ onto $B,$ where $B$ is the
Banach space
of holomorphic functions $\psi(z)$ defined in the lower half-plane
${\mathbb H}^*$ for which
$$||\psi||_B = \sup_{z \in {\mathbb H}^*} |\psi(z)y^2| < \infty.$$
 \end{thm}
\begin{proof}
Note that
$$\alpha(\mu) = {\rm \ Re \ } \left(-\frac{2}{\pi}
\iint_{\mathbb H} \frac{\mu(\zeta)}{\zeta(\zeta-1)(\zeta-x)} d \xi
d
\eta\right),$$
and therefore the condition that $\alpha(\mu)(x)=0$ for all $x$ in
${\mathbb R}$ implies
$$\iint_{\mathbb H} \frac{\mu(\zeta)}{\zeta(\zeta-1)(\zeta-z)} d
\xi d
\eta = 0$$
for all $z$ in the lower half plane.  On taking the third
derivative
with respect to $z,$ we find that
$$\iint_{\mathbb H} \frac{\mu(\zeta)}{(\zeta-z)^4} d \xi d \eta =
0,$$
for all $z$ in the lower half-plane.
Since finite linear combinations of the form
$$ \varphi(\zeta) = \sum_j c_j \frac{1}{(\zeta - z_j)^4},$$
where $z_j$ are points in the lower half-plane,
are dense
in the space of integrable holomorphic quadratic differentials in
the
upper half plane (see \cite{A1}, \cite{Bers1}), we see that
$\alpha(\mu)(x)=0$ for all $x$ in
${\mathbb R}$ implies $\iint_{\mathbb H} \varphi(\zeta) \mu(\zeta)
d
\xi d \eta = 0$ for all $\varphi,$ which implies $\mu$ is in $N.$

Conversely, if $\mu$ is orthogonal to every $\varphi,$ then
$$\iint_{\mathbb H} \frac{\mu(\zeta)}{(\zeta-z)^4} d \xi d \eta =
0,$$
for every $z$ in the lower half plane and, by integrating three
times
and normalizing so that $V_{\mu}(z)$ vanishes at $0,1$ and
$\infty,$ we find that
$$V_{\mu}(x) = - \frac{1}{\pi} \iint_{\mathbb C}
\frac{\hat{\mu}(\zeta)}
{\zeta(\zeta-1)(\zeta-x)} d \xi d \eta = 0,$$
for all $x$ in ${\mathbb R}.$

To see that $\alpha$ is surjective we apply the extension formula
(\ref{BA1}) to the vector field $V(x) \frac{\partial}{\partial x}.$
That is, for given
 $V(x) \frac{\partial}{\partial x}$ representing an element in
${\cal
Z},$ we put
$V(z)=W_1(z)+i W_2(z),$ where
\begin{equation}
\label{BA3}
W_1(x+iy)= \frac{1}{2y} \int_{x-y}^{x+y} V(t)dt,
\end{equation}
and
\begin{equation}
\label{BA4}
W_2(x+iy)= \frac{1}{y} \left\{\int_x^{x+y} V(t)dt -
\int_{x-y}^y h(t) dt \right\}.
\end{equation}
Then it is a routine calculation (see \cite{GS}) to show that
$||\frac{\partial }{\partial \overline{z}} V(z)||_{\infty}<\infty.$

We leave it to the reader to show that the Bers map $\beta$ is an
isomorphism,
and
in particular that
$\beta(\mu)=0$ if, and only if, $V_{\mu}(x)=0$ for every $x$ in
${\mathbb R}.$  A detailed proof may be found in \cite{GL}, p.134,
or in \cite{L},
pp. 111-114.
 \end{proof}
\end{subsection}

\begin{subsection}{Hilbert Transform and Almost Complex Structure}
For a smooth real-valued function $f(x)$ defined on the real axis
with compact support, the Hilbert
transform $Jf$ is normally defined as the principal part of a
divergent
integral.  That is,
\begin{equation}
\label{HT}
\left(Jf\right)(x) = - \frac{1}{\pi} \lim_{\epsilon \rightarrow 0}
\left\{ \int_{-\infty}^{x-\epsilon} f(t) dt +
\int_{x+\epsilon}^{\infty}
f(t) dt \right\}.
\end{equation}
This formula hides a description of the transform in terms of
harmonic conjugates which is invariant under conformal changes of
coordinate. This description has three steps.  The first step, if
it
is possible, is to form the unique harmonic extension
$\tilde{f}(z)$
to the
upper half-plane, characterized by the properties that
$\tilde{f}(z)$ is
harmonic and $\tilde{f}(x)$
coincides with $f(x)$ for $x$ real. Then one forms $\tilde{g}(z),$
which
is unique
up to an additive constant and such that
$\tilde{f}(z)+i\tilde{g}(z)$ is
holomorphic in the upper half-plane.  Finally, $Jf(x)$ is defined
to be
the restriction of $\tilde{g}(z)$ to the real axis.

Of course, the definition of $Jf$ in given this way is determined
only
up to additive constant, but in order for $J$ to be well-defined
on
${\cal Z},$ we need only to define
$JV$ up to the addition of a quadratic polynomial $p(z)=a z^2 + b
z +c.$

We wish to give an alternative  description of the Hilbert
transform on the space
${\cal Z}.$ Since $\alpha$ is surjective, we can assume $V$ is of
the form $$V(x) = - \frac{1}{\pi} \iint_{\mathbb C}
\frac{\mu(\zeta)}
{\zeta(\zeta-1)(\zeta-x)} d \xi d \eta,$$
where $\mu(z)=\overline{\mu(\overline{z})}.$
We shall say that $\mu$ satisfying this equation is {\em
symmetric}.
For symmetric $\mu,$ we define $\hat{\mu}$ to be the
Beltrami coefficient given by the formula
\begin{equation}
\label{symbeltr}
\hat{\mu}(\zeta) =
\left\{ \begin{array}{ll}
\ i \mu(\zeta) & \mbox{ for $\zeta$ in ${\mathbb H}$ } \\
\ & \ \\
-i \mu(\zeta) & \mbox{ for $\zeta$ in ${\mathbb H}^*.$}
\end{array}
\right.
\end{equation}
Then $\hat{\mu}$ is also symmetric and
$$\left(V_{\mu} + i V_{\hat{\mu}}\right)(z)
= \frac{-2}{\pi} \iint_{{\mathbb H}^*}
\frac{\mu(\zeta)}
{\zeta(\zeta-1)(\zeta-z)} d \xi d \eta,$$
where the integration is over the lower half-plane ${\mathbb H}^*.$
It is obvious that this function is holomorphic in the upper
half-plane;
therefore, up to the addition of a quadratic polynomial, $JV_{\mu}$
is the restriction to the real axis of $V_{\hat{\mu}}(z).$

Note that $||\hat{\mu}||_{\infty}=||\mu||_{\infty}.$
This reformulation shows that ${\cal Z}$ is invariant under
$J,$ and in
fact $J$ is an isometry for the infinitesimal \te\ norm on the
tangent
space to \te\ space.
The argument
is
easily modified to show that ${\cal Z}_0$ is also invariant under
$J,$ see (\cite{GS}).  It also shows that the Hilbert transform
applied
to the vector
field $V(x) \frac{\partial}{\partial x}$ corresponds to the mapping
$\mu
\mapsto i \mu$ for Beltrami coefficients given in the upper
half-plane.
Since multiplication by $i$ on Beltrami coefficients determines
the standard almost complex structure
on \te\ space, the Hilbert transform gives the same almost complex
structure.  This observation is
due to Steven Kerckhoff (unpublished, but see \cite{Wolpert}).
\end{subsection}

\begin{subsection}{Complex Structures on Quasi-Fuchsian Space}
The view of the almost complex structures summarized in the
previous
section has been exploited by Giannis Platis \cite{GP} to yield
{\em three} inter-related but distinct
almost complex structures on the quasi-Fuchsian spaces
$QF=QF(\Gamma)$.
These
are complex deformation spaces, whose points are given by
arbitrary
quasiconformal conjugates of a given (cofinite volume) Fuchsian
group $\Gamma$,
acting discretely on the union of $\mathbb H$ and
${\mathbb H}^*$, the complement of the circle
$\overline{\mathbb R}$ inside the Riemann sphere. Such a
group $H\Gamma H^{-1}$ is known as a {\em quasi-Fuchsian group}.
It operates discretely, but not necessarily symmetrically,
on the complement of the quasicircle
$H(\overline{\mathbb R})$; our earlier definition
of the \te\ spaces implies that
$QF(\Gamma)\supset T(\Gamma)$ as a diagonal subset, corresponding
to
q-c conjugates where the mapping
$H$ is given by a symmetric Beltrami coefficient.
Together with a certain hermitian 2-form $\Omega$ defined on the
space
$QF(\Gamma)$, the three anti-involutions
determine a hyper-Kahlerian structure.
One views the tangent space to $QF(\Gamma)$ as a space of (complex)
vector fields $V(x)\frac{\partial}{\partial x}$ which can be
expressed
in the form
$$V(x) = -\frac{1}{\pi} \iint_{\mathbb C}
\frac{\mu(\zeta)}{\zeta(\zeta-1)(\zeta-z)} d \xi d \eta.$$
In this formula, $V(x)$ is usually complex-valued because there is
no assumption about symmetry for $\mu.$

We define $I$ to be the map of vector fields
induced by $\mu \mapsto i \mu.$  Writing
$$\mu = \left\{
\begin{array}{ll}
\mu_1(\zeta) & \mbox{ for $\zeta$ in ${\mathbb H}$ and} \\
\ & \ \\
\mu_2(\zeta) & \mbox{ for $\zeta$ in ${\mathbb H}^*,$}
\end{array}
\right.
$$
we define $V \mapsto J(V)$ to be the map induced by
$$\mu \mapsto  \left\{
\begin{array}{ll}
\overline{\mu_2(\overline{\zeta})} & \mbox{ for $\zeta$ in
${\mathbb H}$ and} \\
\ & \ \\
- \overline{\mu_1(\overline{\zeta})} & \mbox{ for $\zeta$ in
${\mathbb
H}^*.$}
\end{array}
\right.
$$
A simple calculation shows that $IJ=-JI,$ so that if we write
$K=I\circ J,$ then $I^2,$ $J^2$ and $K^2$ are all equal to minus
the identity
and
$IJ=K, JK=I$ and $KI=J.$  Moreover, $K$ restricted to symmetric
Beltrami
coefficients coincides with the almost complex structure defined
(via the Hilbert transform) in the preceding section
on \te\ space.

In \cite{GP}, Platis shows that for finite co-volume Fuchsian
groups, $I, J,$ and $K$ together with the hermitian form
$\Omega$ yield a hyper-Kahlerian structure on $QF(\Gamma).$
The form $\Omega$
is constructed using derivatives of a finite spanning
set of complex length functions
(see for instance \cite{Kou}); it is compatible with the almost
complex structure induced on $QF$ as a product space
$T(\Gamma)\times \overline{T(\Gamma)}$
by the anti-involution $J$,
satisfies a complex analogue of Wolpert's reciprocity formula for
hyperbolic length functions on \te\ space, and restricts
on the diagonal subspace to give the Weil-Petersson
metric on \te\ space. \end{subsection}

\begin{subsection}{Automorphisms are Geometric}
To close this article we return to the rigidity theorem, Theorem
6,
formulated in section 1.8. This result is the analogue for
universal \te\ space of the classical result of H. Royden
\cite{Royden}
and of Earle and Kra \cite{EK2}
that says that any automorphism of the \te\ space of a surface of
genus
greater than $3$ and possibly with a finite number of punctures is
induced by an element of the mapping class
group. That the parallel result holds for any surface of finite
genus
with
a finite number of holes removed was proved in \cite{EG} and for
any
open surface of finite genus by Lakic in \cite{L3}.

Suppose we are given an almost complex diffeomorphism $F$ of
universal
\te\ space, $T.$
Since Kobayashi's metric coincides with \tes\ metric on $T$
\cite{G4},
the
automorphism is an isometry in \tes\ metric, and since \tes\ metric
is
the integral of its infinitesimal form \cite{O}, this means that
if $F([0])=\tau,$ then $F'=dF$ defines an isometry from the tangent
space at $[0]$ to the tangent space at $\tau.$
Since we may select a geometric isomorphism $\rho_h$ such that
$\rho_h
\circ F ([0])=[0]$ and since geometric isomorphisms are isometries,
we
obtain an automorphism $\rho_h \circ F$ which preserves the
basepoint
$[0]$ and induces an isometry on the tangent space ${\cal Z}$ at
$[0]$ to
\te\ space.  One shows that this isometry is necessarily equal to
the
identity
and thus $F=\rho_{h^{-1}},$ that is, every automorphism of \te\
space
is
induced by the action of a quasisymmetric mapping on the boundary
of the
hyperbolic plane.

We outline the key steps in the proof.
One first shows that any isometry $I$ of
${\cal Z}$ with \tes\ infinitesimal metric must be induced by an
isometry
of the predual space $A.$  This result follows from the results of
section 2.3, and in particular from Theorem \ref{thm11}.  Such an
isometry must preserve the closed subspace ${\cal Z}_0$ and
therefore
it is equal to the second dual of its restriction to ${\cal Z}_0.$
Thus
there is an
isometry $\hat{I}$ of ${\cal A}$
such that $I$ is the dual of $\hat{I}$ under the natural pairing
between
$A$ and ${\cal Z}.$
Then one shows that $\hat{I}$ is induced by the composition of
multiplication
by a complex constant $c$ (of modulus 1 since it is an isometry)
with a conformal map.
That is,
$$\hat{I}(\varphi) = c \varphi(f(z))f'(z)^2,$$
where $f$ is a conformal self-map of the base Riemann surface.  For
universal \te\ space, the base Riemann surface is the upper
half-plane and so $f$ is a real M\"obius transformation in this case.
Finally, one shows that the constant $c$ is equal to
$1:$  for this step, see
\cite{EG} or \cite{GL}.
\end{subsection}

\end{section}

\bibliography{newver}

\def\noopsort#1{}
\def\printfirst#1#2{#1}
\def\singleletter#1{#1}
\def\switchargs#1#2{#2#1}
\def\bibsameauth{\leavevmode\vrule height .1ex
  depth 0pt width 2.3em\relax\,}

\end{document}